\documentclass[11pt]{amsart}

\usepackage{epsfig,amssymb}

\makeatletter
\def\eqalign#1{\null\vcenter{\def\\{\cr}\openup\jot\m@th
  \ialign{\strut$\displaystyle{##}$\hfil&$\displaystyle{{}##}$\hfil
      \crcr#1\crcr}}\,}
\makeatother

\newcommand{\be}{\begin{equation}} 
\newcommand{\ee}{\end{equation}}
\newcommand{\beq}{\begin{eqnarray}}
\newcommand{\eeq}{\end{eqnarray}}
\newcommand{\bt}{\beta}
\newcommand{\bl}{\begin{lemma}}
\newcommand{\el}{\end{lemma}}
\newcommand{\bm}{\begin{pmatrix}}
\renewcommand{\em}{\end{pmatrix}}
\newcommand{\bml}{\begin{multline}}
\newcommand{\eml}{\end{multline}}
\newcommand{\ba}{\begin{array}}
\newcommand{\ea}{\end{array}}

\newcommand{\la}{\label}
\newcommand{\ci}{\cite}

\newcommand{\de}{\delta}
\newcommand{\De}{\Delta}
\newcommand{\al}{\alpha}
\newcommand{\ga}{\gamma}
\newcommand{\Ga}{\Gamma}

\newcommand{\si}{\sigma}
\newcommand{\Si}{\Sigma}
\newcommand{\om}{\omega}
\newcommand{\Om}{\Omega}
\newcommand{\lb}{\lambda}
\newcommand{\ze}{\zeta}

\renewcommand{\th}{\theta}

\newcommand{\ep}{\varepsilon }
\newcommand{\wt}{\widetilde}
\newcommand{\wh}{\widehat}

\newcommand{\bi}{\bibitem}

\newfont{\msbm}{msbm10 scaled\magstep1}%blackboardbold
\newfont{\msbms}{msbm7 scaled\magstep1} %blackboardbold   subscript

\newcommand{\bbr}{\mbox{$\mbox{\msbm R}$}}

\newcommand{\bbc}{\mbox{$\mbox{\msbm C}$}}

\newcommand{\bbz}{\mbox{$\mbox{\msbm Z}$}}

\newtheorem{theorem}{Theorem}[section]
\newtheorem{lemma}{Lemma}[section]
\newtheorem{proposition}[theorem]{Proposition}

\theoremstyle{definition}

\theoremstyle{remark}
\newtheorem{remark}[theorem]{Remark}

\numberwithin{equation}{section}

\begin{document}
%\def\wt{\widetilde}
%\hfill {\small February 28, 2008}
\title[Toeplitz determinants II]{On the asymptotics of a Toeplitz determinant with singularities}
\author{P. Deift}
\address{Courant Institute of Mathematical Sciences, New York, NY 10003, USA}

\author{A. Its}
\address{Department of Mathematical Sciences,
Indiana University -- Purdue University  Indianapolis,
Indianapolis, IN 46202-3216, USA}
%\thanks{The second author was supported in part by NSF Grant \#DMS-0401009.}
%
\author{I. Krasovsky}
\address{Department of Mathematics,
Imperial College,
London SW7 2AZ, United Kingdom}
%
%\thanks{The third author was supported in part by EPSRC Grant EP/E022928/1.}

%\date{}

\begin{abstract}
We provide an alternative proof of the classical single-term asymptotics for Toeplitz determinants
whose symbols possess Fisher-Hartwig singularities. We also relax the smoothness conditions on 
the regular part of the symbols and obtain an estimate for the error term in the asymptotics. 
Our proof is based on the Riemann-Hilbert analysis of the related 
systems of orthogonal polynomials and on differential identities for Toeplitz determinants. 
The result discussed in this paper is crucial for the proof of the asymptotics 
in the general case of Fisher-Hartwig singularities and extensions to Hankel and Toeplitz+Hankel 
determinants in \ci{DIKt1}. 
\end{abstract}

\maketitle

\section{Introduction}
Let $f(z)$ be a complex-valued function integrable over the unit circle. Denote its Fourier coefficients
\[
f_j={1\over 2\pi}\int_0^{2\pi}f(e^{i\theta})e^{-i j\theta}d\theta,\qquad j=0,\pm1,\pm2,\dots
\]
We are interested in the $n$-dimensional Toeplitz determinant with symbol $f(z)$,
\be\la{TD}
D_n(f(z))=\det(f_{j-k})_{j,k=0}^{n-1},\qquad n\ge 1,
\ee
where $f(e^{i\theta})$ has a fixed number of 
Fisher-Hartwig singularities \ci{FH,L}, i.e., 
$f$ has the following form on the unit circle:
\be\la{fFH}
f(z)=e^{V(z)} z^{\sum_{j=0}^m \bt_j} 
\prod_{j=0}^m  |z-z_j|^{2\al_j}g_{z_j,\bt_j}(z)z_j^{-\bt_j},\qquad z=e^{i\th},\qquad
\theta\in[0,2\pi),
\ee
for some $m=0,1,\dots$,
where
\begin{eqnarray}
&z_j=e^{i\th_j},\quad j=0,\dots,m,\qquad
0=\th_0<\th_1<\cdots<\th_m<2\pi;&\la{z}\\
&g_{z_j,\bt_j}(z)\equiv g_{\bt_j}(z)=
\begin{cases}
e^{i\pi\bt_j}& 0\le\arg z<\th_j\cr
e^{-i\pi\bt_j}& \th_j\le\arg z<2\pi
\end{cases},&\la{g}\\
&\Re\al_j>-1/2,\quad \bt_j\in\bbc,\quad j=0,\dots,m,&
\end{eqnarray}
and $V(e^{i\theta})$ is a sufficiently smooth function on the unit circle (see below). 
The condition on the $\al_j$'s insures integrability.
Note that a single Fisher-Hartwig singularity at $z_j$ consists of 
a root-type singularity
\be\la{za}
|z-z_j|^{2\al_j}=\left|2\sin\frac{\th-\th_j}{2}\right|^{2\al_j}
\ee
and a jump $e^{i\pi\bt}\to e^{-i\pi\bt}$.
We  assume that $z_j$, $j=1,\dots,m$, are genuine singular points, i.e.,
either $\al_j\neq 0$ or $\bt_j\neq 0$. However, we always include $z_0=1$ explicitly 
in (\ref{fFH}),
even when $\al_0=\bt_0=0$: this convention was adopted in \cite{DIKt1} in order to 
facilitate the application 
of our Toeplitz methods to Hankel determinants.
Note that $g_{\bt_0}(z)=e^{-i\pi\bt_0}$.
Observe that for each $j\neq 0$, $z^{\beta_j} g_{\beta_j}(z)$ 
is continuous at $z=1$, and so for each $j$ 
each ``beta'' singularity produces a jump only at the point $z_j$.
The factors $z_j^{-\bt_j}$ are singled out to simplify comparisons 
with the existing literature. Indeed, (\ref{fFH}) with the notation 
$b(\th)=e^{V(e^{i\theta})}$ is exactly the symbol considered 
in \ci{FH,B,B2,BS,BS2,BS3,ES,Ehr,W}. 
However, we write the symbol in a form with
$z^{\sum_{j=0}^m \bt_j}$ factored out. The representation (\ref{fFH})
is more natural for our analysis. 

On the unit circle, $V(z)$ is represented by its Fourier expansion:
\be\la{fourier}
V(z)=\sum_{k=-\infty}^\infty V_k z^k,\qquad 
V_k={1\over 2\pi}\int_0^{2\pi}V(e^{i\th})e^{-ki\th}d\th.
\ee
The canonical Wiener-Hopf factorization of $e^{V(z)}$ is given by 
\be\la{WienH}
e^{V(z)}=b_+(z) e^{V_0} b_-(z),\qquad b_+(z)=e^{\sum_{k=1}^\infty V_k z^k},
\qquad b_-(z)=e^{\sum_{k=-\infty}^{-1} V_k z^k}.
\ee

Define a seminorm
\be\la{norm}
|||\bt|||=\max_{j,k}|\Re\bt_j-\Re\bt_k|,
\ee
where the indices $j,k=0$ are omitted if $z=1$ is not a singular point, i.e. if $\al_0=\bt_0=0$.
If $m=0$, set $|||\bt|||= 0$.

In this paper we consider the asymptotics of $D_n(f)$, $n\to\infty$, 
in the case $|||\bt|||<1$. The case $|||\bt|||\ge 1$
was addressed in \ci{DIKt1}.
The asymptotic behavior of $D_n(f)$ has been studied by many authors (see \ci{DIKt1,Ehr} 
for a review). 
An expansion of $D_n(f)$ for the case $V\in C^\infty$,  $|||\bt|||<1$, was obtained by Ehrhardt in \ci{Ehr}.
The aim of this paper is to provide an alternative proof of this result based on differential 
identities for 
$D_n(f)$  and a Riemann-Hilbert-Problem analysis of the corresponding system of orthogonal 
polynomials. We also obtain estimates for the error term and extend the validity
of the result to less smooth $V$. 
The analysis of  $D_n(f)$ for $|||\bt|||<1$ plays a crucial role in the analysis of $D_n(f)$ 
for $|||\bt|||\ge 1$ (see \ci{DIKt1}).

We prove
\begin{theorem}\la{asTop}
Let $f(e^{i\theta})$ be defined in (\ref{fFH}), $|||\bt|||<1$, $\Re\al_j>-1/2$,  
$\al_j\pm\bt_j\neq -1,-2,\dots$ for $j,k=0,1,\dots,m$, and 
let $V(z)$ satisfy the condition (\ref{Vcond}), (\ref{s-main0}) below. 
Then as $n\to\infty$,
\begin{multline}\la{asD}
D_n(f)=\exp\left[nV_0+\sum_{k=1}^\infty k V_k V_{-k}\right]
\prod_{j=0}^m b_+(z_j)^{-\al_j+\bt_j}b_-(z_j)^{-\al_j-\bt_j}\\
\times
n^{\sum_{j=0}^m(\al_j^2-\bt_j^2)}\prod_{0\le j<k\le m}
|z_j-z_k|^{2(\bt_j\bt_k-\al_j\al_k)}\left({z_k\over z_j e^{i\pi}}
\right)^{\al_j\bt_k-\al_k\bt_j}\\
\times
\prod_{j=0}^m\frac{G(1+\al_j+\bt_j) G(1+\al_j-\bt_j)}{G(1+2\al_j)}
\left(1+o(1)\right),
\end{multline}
where
$G(x)$ is Barnes' $G$-function. The double product over $j<k$ is set to $1$
if $m=0$.
\end{theorem}

\begin{remark} As indicated above, this result was first obtained by Ehrhardt in the case $V\in C^{\infty}$.
We prove the theorem for $V(z)$ satisfying the smoothness condition
\be\la{Vcond}
\sum_{k=-\infty}^\infty |k|^s |V_k|<\infty,
\ee
where
\be\la{s-main0}
s>
\frac{1+\sum_{j=0}^m\left[(\Im\al_j)^2+(\Re\bt_j)^2\right]}{1-|||\bt|||}.
\ee
\end{remark}

\begin{remark}
In the case of a single singularity,
i.e., when $m=1$ and $\al_0=\bt_0=0$, or when $m=0$,
the seminorm $|||\bt|||=0$, and
the theorem implies that the asymptotic form (\ref{asD}) holds for all
\be
\Re\al_m>-{1\over 2},\qquad \bt_m\in\bbc,\qquad \al_m\pm\bt_m\neq -1,-2,\dots
\ee
In fact, if there is only one singularity, say $m=0$, 
and $V\equiv 0$, an explicit formula was found by B\"ottcher and Silbermann \ci{BS}
for $D_n(f)$ for any $n$ in terms of the G-functions:
\be\label{eqBS}
\begin{aligned}
D_n(f)=&\frac{G(1+\al_0+\beta_0)G(1+\al_0-\beta_0)}{G(1+2\al_0)}
\frac{G(n+1)G(n+1+2\al_0)}{G(n+1+\al_0+\beta_0)G(n+1+\al_0-\beta_0)},\\
&\Re\al_0>-{1\over 2},\qquad\al_0\pm\beta_0\neq -1,-2,\dots,\qquad n\ge 1,
\end{aligned}
\ee
and (\ref{asD}) can then be read off from the known asymptotics of the G-function (see, e.g., \cite{Barnes}). 
\end{remark}

\begin{remark}\la{error}
Assume that the function $V(z)$ is sufficiently smooth, i.e. 
such that $s$ in (\ref{Vcond}) is, in addition to satisfying (\ref{s-main0}), 
sufficiently large in comparison with $|||\bt|||$.
Then we show that the error term $o(1)=O(n^{|||\bt|||-1})$ in (\ref{asD}). 
In particular, the error term $o(1)=O(n^{|||\bt|||-1})$ if $V(z)$ is analytic in a neighborhood
of the unit circle. 
Moreover, for analytic $V(z)$, our methods would allow us to calculate, in principle, 
the full asymptotic expansion rather than just the leading term presented in 
(\ref{asD}). Various regularity properties of the expansion (uniformity, differentiability) in
compact sets of parameters satisfying $\Re\al_j>-1/2$, $|||\bt|||<1$, $\al_j\pm\bt_j\neq -1,-2,\dots,
\th_j\neq \th_k$, are easy to deduce from our analysis.
\end{remark}

%\begin{remark}
%If all $\Re\bt_j\in(-1/2,1/2]$ or all $\Re\bt_j\in[-1/2,1/2)$,
%the conditions $\al_j\pm\bt_j\neq -1,-2,\dots$ are satisfied automatically
%as $\Re\al_j>-1/2$.
%\end{remark}

\begin{remark}\la{degen1}
Since $G(-k)=0$, $k=0,1,\dots$, the formula (\ref{asD}) no longer 
represents the leading asymptotics
if $\al_j+\bt_j$ or $\al_j-\bt_j$ is a negative integer for some $j$.
Although our method applies,
we do not address these cases in this paper. It is simply a matter of going deeper in the asymptotic
expansion for $D_n(f)$. It can happen that $D_n(f)$ vanishes to all orders (cf. discussion of the Ising model
at temperatures above the critical temperature in \ci{DIKrev}), and $e^{-nV_0}D_n(f)$ is exponentially
decreasing.
\end{remark}

We prove Theorem \ref{asTop} in the following way. We begin by deriving differential identities
for  the logarithm of $D_n(f)$ in Section \ref{diffid}, in the spirit 
of \ci{DeiftInOp,Ksine,Khankel,IK,DIKZ,DIKAiry}, 
utilizing the polynomials orthogonal with respect to the weight $f(z)$ on the unit circle. 
Then, assuming that $V(z)$ is analytic in a neighborhood of the unit circle, 
we analyze in Section \ref{RHa} 
the asymptotics of these polynomials using Riemann-Hilbert/steepest-descent methods as in \ci{DIKt1}.
This gives in turn the asymptotics of the differential identities from which 
the formula (\ref{asD}) follows in the $V\equiv 0$ case  by integration w.r.t. $\al_j$, $\bt_j$
in Section \ref{Toep51}. 
However,  the error term that results is of order $n^{2|||\bt|||-1}\ln n$ (see (\ref{542})),
which is asymptotically small only for $|||\bt|||<1/2$, rather than in the full 
range $|||\bt|||<1$. 
To prove (\ref{asD}) for all $|||\bt|||<1$, we need a finer analysis of cancellations 
in the Riemann-Hilbert problem as $n\to\infty$.
We carry this out in Section \ref{pFHext} and reduce the leading order terms in $D_n(f)$ to 
a telescopic form (see (\ref{newestD})),
which leads to a uniform bound on $D_n(f)n^{\sum_{j=0}^m(\bt_j^2-\al_j^2)}$
for large $n$ which is valid for all $|||\bt|||<1$ and away from the points 
$\al_j\pm\bt_j=-1,-2,\dots$.
We then apply Vitali's theorem together with the previous result for
$|||\bt|||<1/2$. This proves Theorem \ref{asTop} in the $V\equiv 0$ case as desired, with the
error term of order $n^{|||\bt|||-1}$.
In Sections \ref{secV}, \ref{secVGenAn}, we then extend the result to the case of analytic $V\not\equiv 0$
by applying another differential identity from Section \ref{diffid}. 

Ehrhardt proves (\ref{asD}) using a ``localization'' or ``separation'' technique, introduced
by Basor in \ci{B2}, in which the effect of adding in Fisher-Hartwig singularities one at a time,
is controlled. One may also view our approach as a ``separation'' technique, but in contrast to \ci{B2,Ehr},
we add in the Fisher-Hartwig singularities, as well as the regular term $e^{V(z)}$, in a continuous fashion.

\begin{remark}
An alternative approach to proving Theorem \ref{asTop} in the $V$-analytic case is to 
apply, ab initio, the finer analysis of Section \ref{pFHext} to the orthogonal polynomials which 
appear in the differential identities.  
This approach is more direct, but is considerably more involved technically.
The analysis is resolved, as above, by reduction of the problem to an appropriate telescopic form.
\end{remark}

Finally in Section \ref{secVext}, we extend our result to the case 
when $V(z)$ is not analytic and only satisfies 
the smoothness condition (\ref{Vcond}), (\ref{s-main0}). 
We approximate such $V(z)$ by trigonometric polynomials $V^{(n)}(z)=\sum_{k=-p}^p V_k z^k$
with an appropriate $p=p(n)$ and modify the Riemann-Hilbert analysis accordingly. 
This produces  asymptotics of a Toeplitz
determinant in whose symbol $f^{(n)}(z)$ the function $V(z)$ is replaced by $V^{(n)}(z)$.
We then use the Heine representation of Toeplitz determinants by multiple integrals 
to show that $D_n(f^{(n)})$ approximates $D_n(f)$ as $n\to\infty$, sufficiently strongly 
to conclude Theorem \ref{asTop} in the general case.

\section{Riemann-Hilbert problem}\la{RHsection}
In this section we formulate a Riemann-Hilbert problem (RHP) 
for the polynomials orthogonal on the unit circle (which, oriented in the positive direction, 
we denote $C$) 
with weight $f(z)$ given by (\ref{fFH}).
We use this RHP in Section \ref{RHa} to find the asymptotics of the polynomials 
in the case of analytic $V(z)$.
Suppose that all $D_k(f)\neq 0$, $k=k_0,k_0+1\dots$,
for some sufficiently large $k_0$ (see discussion below).
Then the polynomials 
$\phi_k(z)=\chi_k z^k+\cdots$, $\widehat\phi_k(z)=\chi_k z^{k}+\cdots$ 
of degree $k$, $k=k_0,k_0+1,\dots$, satisfying the orthogonality conditions 
\begin{multline}\la{or0}
{1\over 2\pi}\int_0^{2\pi}\phi_k(z)z^{-j}f(z)d\theta=\chi_k^{-1}\de_{jk},\qquad
{1\over 2\pi}\int_0^{2\pi}\widehat\phi_k(z^{-1})z^j f(z)d\theta=
\chi_k^{-1}\de_{jk},\\
z=e^{i\theta},\qquad j=0,1,\dots,k,
\end{multline}
exist and are given by the following expressions:
\be\la{ef1}
\phi_k(z)={1\over\sqrt{D_k D_{k+1}}}
\left| 
\begin{matrix}
f_{00}& f_{01}& \cdots & f_{0k}\cr
f_{10}& f_{11}& \cdots & f_{1k}\cr
\vdots & \vdots &  & \vdots \cr
f_{k-1\,0} & f_{k-1\,1} & \cdots & f_{k-1\,k} \cr
1& z& \cdots & z^k
\end{matrix}
\right|,
\ee
\be\la{ef2}
\widehat\phi_k(z^{-1})={1\over\sqrt{D_k D_{k+1}}}
\left| 
\begin{matrix}
f_{00}& f_{01}& \cdots & f_{0\,k-1}& 1\cr
f_{10}& f_{11}& \cdots & f_{1\,k-1}& z^{-1}\cr
\vdots & \vdots &  & \vdots & \vdots\cr
f_{k0} & f_{k1} & \cdots & f_{k\,k-1}& z^{-k}
\end{matrix}
\right|,
\ee
where
\[ 
f_{st}={1\over 2\pi}\int_0^{2\pi}f(z)z^{-(s-t)}d\theta,\quad 
s,t=0,1,\dots,k.
\]
We obviously have
\be\la{chiD}
\chi_k=\sqrt{D_k\over D_{k+1}}.
\ee

Consider the following $2\times 2$ matrix valued function $Y^{(k)}(z)\equiv Y(z)$, $k\ge k_0$:
\begin{equation} \label{RHM}
    Y^{(k)}(z) =
    \begin{pmatrix}
\chi_k^{-1}\phi_k(z) & 
\chi_k^{-1}\int_{C}{\phi_k(\xi)\over \xi-z}
{f(\xi)d\xi \over 2\pi i \xi^k} \cr
-\chi_{k-1}z^{k-1}\wh\phi_{k-1}(z^{-1}) & 
-\chi_{k-1}\int_{C}{\wh\phi_{k-1}(\xi^{-1})\over \xi-z}
{f(\xi)d\xi \over 2\pi i \xi}
    \end{pmatrix}.
\end{equation}

It is easy to verify that 
$Y(z)$ solves the following Riemann-Hilbert problem:
\begin{enumerate}
    \item[(a)]
        $Y(z)$ is  analytic for $z\in\bbc \setminus C$.
    \item[(b)] 
Let $z\in C\setminus\cup_{j=0}^m z_j$.
$Y$ has continuous boundary values
$Y_{+}(z)$ as $z$ approaches the unit circle from
the inside, and $Y_{-}(z)$, from the outside, 
related by the jump condition
\begin{equation}\label{RHPYb}
            Y_+(z) = Y_-(z)
            \begin{pmatrix}
                1 & z^{-k}f(z) \cr
                0 & 1
             \end{pmatrix},
            \qquad\mbox{$z\in C\setminus\cup_{j=0}^m z_j$.}
        \end{equation}
    \item[(c)]
        $Y(z)$ has the following asymptotic behavior at infinity:
        \begin{equation} \label{RHPYc}
            Y(z) = \left(I+ O \left( \frac{1}{z} \right)\right)
            \begin{pmatrix}
                z^{k} & 0 \cr
                0 & z^{-k}\end{pmatrix}, \qquad \mbox{as $z\to\infty$.}
        \end{equation}
\item [(d)] 
As $z\to z_j$, $j=0,1,\dots,m$, $z\in\bbc\setminus C$,
\be\label{RHPYd}
Y(z)=
\begin{pmatrix} O(1) & O(1)+O(|z- z_j|^{2\al_j})\cr 
O(1) & O(1)+O(|z- z_j|^{2\al_j})
\end{pmatrix},
\qquad \mbox{if $\al_j\neq 0$},
\ee
and
\be\label{RHPYe}
Y(z)=
\begin{pmatrix} O(1) & O(\ln|z- z_j|)\cr 
O(1) & O(\ln|z-z_j|)
\end{pmatrix},
\qquad \mbox{if $\al_j=0$, $\bt_j\neq 0$}.
\ee
\end{enumerate}
(Here and below $O(a)$ stands for $O(|a|)$.) If $\al_0=\bt_0=0$, $Y(z)$ is bounded at $z=1$.

A general fact that orthogonal polynomials can be so represented as a solution
of a Riemann-Hilbert problem was noticed in \cite{FIK} for polynomials
on the line and extended to polynomials on the circle in \cite{BDJ}.
This fact is important because it turns out that the RHP can be efficiently analyzed
for large $k$ by a steepest-descent type method found in \ci{DZ} and developed 
further in many subsequent works. Thus, we first 
find the solution to the problem
(a)--(d) for large $k$ (applying this method)
and then interpret it as the asymptotics of 
the orthogonal polynomials by (\ref{RHM}).

Recall the Heine representation for a Toeplitz determinant:
\be\la{ir}
D_{n}(f) = \frac{1}{(2\pi)^{n}n!}
\int_0^{2\pi}\cdots\int_0^{2\pi}
\prod_{1\leq j < k \leq n}|e^{i\phi_{j}} - e^{i\phi_{k}}|^{2}
\prod_{j=1}^{n} f(e^{i\phi_j})d\phi_{j}.
\ee
If $f(z)$ is positive on the unit circle, it follows from (\ref{ir}) that
$D_k(f)>0$, $k=1,2,\dots$, i.e. $k_0=1$. In the general case, let $\Lambda$ be a compact subset
in the subset $|||\bt|||<1$, $\al_j\pm\bt_j\neq -1,-2,\dots$ of the parameter
space 
${\mathcal P}=\{(\al_0,\bt_0,\dots,\al_m,\bt_m):\, \al_j,\bt_j\in\bbc,\, \Re\al_j>-1/2\}$. 
We will show in Section \ref{RHa} that the
RHP (a)--(d) is solvable in $\Lambda$, in particular $\chi_k$ are finite and nonzero,
for all sufficiently large $k$ ($k\ge k_0(\Lambda)$). Let $\Om_{k_0}$ be the set of parameters
in ${\mathcal P}$ such that $D_k(f)=0$ for some $k=1,2,\dots,k_0-1$. We will then have 
$D_k(f)\neq 0$, $k=1,2,\dots$ for all points in $\Lambda\setminus\Om_{k_0}$.
Note that $D_n(f)$ depends analytically on $\al_j$, $\bt_j$ in ${\mathcal P}$:
this is true, in particular, on (the interior of) $\Lambda\setminus\Om_{k_0}$.

The solution to the RHP (a)--(d) is unique. Note first
that $\det Y(z)=1$. Indeed, from the conditions on $Y(z)$, 
$\det Y(z)$ is analytic across the unit circle, has all singularities 
removable, and tends to $1$ as $z\to\infty$. 
It is then identically $1$ by Liouville's theorem. Now  
if there is another solution $\wt Y(z)$, we easily obtain by 
Liouville's theorem that $\wt Y(z)Y(z)^{-1}\equiv 1$.

%%%%%%%%%%%%%%%%%%%% diffid %%%%%%%%%%%%%%%%%%%%%

\section{Differential identities}\la{diffid}
In this section we derive expressions for the derivative 
$(\partial/\partial\ga)\ln D_n(f(z))$, where either $\ga=\al_j$ or $\ga=\bt_j$, $j=0,1,\dots,m$,
in terms of the matrix elements of (\ref{RHM}). These will be exact 
differential identities valid for all $n=1,2,\dots$ (see Proposition \ref{idDab} below), provided 
all the $D_n(f)\neq 0$.
We will use these expressions in Section \ref{Toep51}
to obtain the asymptotics (\ref{asD}) in the case $V\equiv 0$, $|||\bt|||<1/2$
(improved to $|||\bt|||<1$ in Section \ref{pFHext}).
Furthermore, in this section we will
derive a differential identity (see  Proposition \ref{idDt} below) which will enable us 
in Sections \ref{secV}, \ref{secVGenAn} to extend the results to analytic $V\not\equiv 0$
(improved to sufficiently smooth $V$ in Section \ref{secVext}).

Set $D_0\equiv 1$,  $\phi_0(z)\equiv\wh\phi_0(z)\equiv 1$, and
suppose that $D_n(f)\neq 0$ for all $n=1,2,\dots$.
Then the orthogonal polynomials (\ref{ef1}),
(\ref{ef2}) exist and are analytic in the $\al_j$'s, $\bt_j$'s for all $k=1,2,\dots$. 
Moreover, (\ref{chiD}) implies that
\be\la{Dprod}
D_n(f(z))=\prod_{j=0}^{n-1}\chi_j^{-2}.
\ee

Note that by orthogonality,
\begin{multline}\la{chidif1}
\frac{1}{2\pi}\int_0^{2\pi}{\partial\phi_j(z)\over\partial\ga} \wh\phi_j(z^{-1})
f(z)d\th=\\ 
\frac{1}{2\pi}\int_0^{2\pi}\left({\partial\chi_j\over\partial \ga} z^j + 
\mbox{polynomial of degree $j-1$}\right)\wh\phi_j(z^{-1})f(z)d\th=
{1\over\chi_j}{\partial\chi_j\over \partial\ga}.
\end{multline}
Similarly,
\be\la{chidif2}
\frac{1}{2\pi}\int_0^{2\pi}\phi_j(z){\partial\wh\phi_j(z^{-1})\over\partial\ga}
f(z)d\th= 
{1\over\chi_j}{\partial\chi_j\over\partial \ga}.
\ee
Therefore, using equation (\ref{Dprod}), we obtain
\be\la{premid}
{\partial \over\partial \ga}\ln D_n(f(z))=
{\partial \over\partial \ga}\ln 
\prod_{j=0}^{n-1}\chi_j^{-2}=
-2\sum_{j=0}^{n-1}{{\partial\chi_j\over\partial \ga}\over \chi_j}=
-\frac{1}{2\pi}\int_0^{2\pi}
{\partial \over\partial \ga}
\left(\sum_{j=0}^{n-1} \phi_j(z)\wh\phi_j(z^{-1})\right) f(z)d\th.
\ee

Using here the Christoffel-Darboux identity (see, e.g., Lemma 2.3 of \cite{DIKt1})
\[
\sum_{k=0}^{n-1}\widehat\phi_k(z^{-1})\phi_k(z)=
-n\phi_n(z)\widehat\phi_n(z^{-1}) +z\left(\widehat\phi_n(z^{-1}){d\over dz}\phi_n(z)-
\phi_n(z){d\over dz}\widehat\phi_n(z^{-1})\right),
\]
and then orthogonality,
we can write 
\be\la{mid}
{\partial \over\partial \ga}\ln D_n(f(z))=
2n {{\partial\chi_n\over\partial \ga}\over \chi_n}+
\frac{1}{2\pi}\int_0^{2\pi} 
{\partial \over\partial \ga}
\left(\phi_n(z){d\wh\phi_n(z^{-1})\over dz}
-\wh\phi_n(z^{-1}){d\phi_n(z)\over dz}\right)z f(z) d\th.
\ee

Writing out the derivative w.r.t. $\ga$ in the integral and using orthogonality, we obtain:
\be\la{mid2}
{\partial \over\partial \ga}\ln D_n(f(z))=I_1-I_2,
\ee
where 
\be
I_1=\frac{1}{2\pi i}\int_0^{2\pi}
{\partial\phi_n(z) \over\partial \ga}
{\partial\wh\phi_n(z^{-1}) \over\partial \th} f(z) d\th,\qquad
I_2=\frac{1}{2\pi i}\int_0^{2\pi}
{\partial\phi_n(z) \over\partial \th}
{\partial\wh\phi_n(z^{-1}) \over\partial \ga} f(z) d\th.
\ee

It turns out that the particular structure of Fisher-Hartwig singularities
allows us to reduce (\ref{mid2}) to a local formula, i.e. to replace the
integrals by the polynomials (and their Cauchy transforms) 
evaluated only at several points (cf. \ci{IK,Khankel}).

Let us encircle each of the points
$z_{j}$ by a sufficiently small disc,
\be\la{U}
U_{z_{j}} = \left\{z:|z-z_{j}| < \ep\right\}.
\ee
Denote
\be\la{Cep}
C_{\ep} = \cup_{j=0}^{m}\left(U_{z_{j}}\cap C\right).
\ee

We now integrate $I_1$ by parts. First assume that $V(z)\equiv 0$.
Then, using the expression 
\[
{\partial f(z) \over\partial \th}=
\sum_{j=0}^m \left(\al_j 
\mbox{cot\;}\frac{\th-\th_j}{2}+i\bt_j\right)f(z)=
\left(\sum_{j=0}^m \al_j \frac{z+z_j}{z-z_j}+\bt_j\right)i f(z),
\]
we obtain:
\begin{multline}\la{I1}
I_1=
 -\chi_n^{-1}{\partial\chi_n\over\partial \ga}
\left(n+\sum_{j=0}^m\bt_j\right)-
\lim_{\ep\to 0}\left[\frac{1}{2\pi}\int_{C\setminus C_\ep}
{\partial\phi_n(z) \over\partial \ga}
\wh\phi_n(z^{-1}) \left(
\sum_{j=0}^m\al_j \frac{z+z_j}{z-z_j}\right)
f(z)d\th\right.\\
\left.
-\frac{1}{2\pi i}
\sum_{j=0}^m {\partial\phi_n(z_j) \over\partial \ga}\wh\phi_n(z_j^{-1})
(f(z_je^{-i\ep})-f(z_je^{i\ep}))\right],
\end{multline}
where the integration is over $C\setminus C_{\ep}$ in the positive direction around
the unit circle.

Note that by adding and subtracting $ {\partial\phi_n(z_j) \over\partial \ga}$,
and by using orthogonality
we can write
\begin{multline}
\int_{C\setminus C_\ep}
{\partial\phi_n(z) \over\partial \ga}
\wh\phi_n(z^{-1}) \frac{z+z_j}{z-z_j}
f(z)d\th=\\
\int_{C\setminus C_\ep}\wh\phi_n(z^{-1}){{\partial\phi_n(z) \over\partial \ga}-
{\partial\phi_n(z_j) \over\partial \ga}\over z-z_j}
(z+z_j) f(z)d\th+
{\partial\phi_n(z_j) \over\partial \ga}
 \int_{C\setminus C_\ep}\wh\phi_n(z^{-1}){2z_j\over z-z_j}  f(z)d\th+ O(\ep^{2\Re\al_j+1}).
\end{multline}
Obviously, the fraction in the first integral on the r.h.s. is 
a polynomial in $z$ of degree $n-1$ with leading coefficient 
$\partial\chi_n/\partial\ga$.
Therefore, the integral equals $2\pi{\partial\chi_n \over\partial \ga}/\chi_n$
up to $O(\ep^{2\Re\al_j+1})$.
The second integral can be written in terms of the element $Y_{22}$ of (\ref{RHM})
for $z\to z_j$. Let us estimate therefore the following expression for $\al_j\neq 0$:
\begin{multline}\la{intdiff}
 \int_C\wh\phi_n(s^{-1}){2z_j f(s)\over s-z}\frac{ds}{is}-
 \lim_{\ep\to 0}\left[
\int_{C\setminus C_\ep}\wh\phi_n(s^{-1}){2z_j f(s)\over s-z_j}\frac{ds}{is}-
{1\over i\al_j}\wh\phi_n(z_j^{-1})(f(z_je^{-i\ep})-f(z_je^{i\ep})\right],
\\
z\to z_j,\qquad |z|>1.
\end{multline}
This difference tends to zero as $z\to z_j$,  
for $\Re\al_j>0$. When $\Re\al_j<0$, it is a growing function as $z\to z_j$,  
and when $\Re\al_j=0$, $\Im\al_j\neq 0$, an oscillating one.
The analysis is similar to that of Section 3 in \ci{Khankel}. For future use, we now
fix an analytical continuation of the absolute value, namely, write
for $z$ on the unit circle,
\be\la{abszzj}
|z-z_j|^{\al_j}=(z-z_j)^{\al_j/2}(z^{-1}-z_j^{-1})^{\al_j/2}=
\frac{(z-z_j)^{\al_j}}{(z z_j e^{i\ell_j})^{\al_j/2}},\qquad z=e^{i\th},
\ee
where $\ell_j$ is found from the condition that the argument of the above function
is zero on the unit circle.
Let us fix the cut of $(z-z_j)^{\al_j}$ going along the line $\th=\th_j$ from $z_j$
to infinity. Fix the branch by the condition that on the line going from $z_j$ to 
the right parallel to the real axis, $\arg (z-z_j)=2\pi $. For $z^{\al_j/2}$ in the denominator, 
$0<\arg z<2\pi$. If $z_0=1$ let $0<\arg (z-1)<2\pi$.
(This choice will
enable us to use the standard asymptotics for a confluent hypergeometric function
in the RH analysis in Section \ref{RHa} below.) 
Then, a simple consideration of triangles shows that 
\be
\ell_j=\begin{cases}
3\pi,& 0<\th<\th_j\cr
\pi, & \th_j<\th<2\pi
\end{cases}
\ee
Thus (\ref{abszzj}) is continued analytically to neighborhoods of the arcs $0<\th<\th_j$, and
$\th_j<\th< 2\pi$.
We now analyze (\ref{intdiff}) in the same way that equation (26) was analyzed in \ci{Khankel}.
For this analysis, however, we will need two other choices
of the function $ (z-z_j)^{2\al_j}$: one choice with the cut going a short distance clockwise along the unit circle $C$ from $z_j$, and another, with the cut going a short distance anticlockwise along $C$
from $z_j$. Let $c_j$ and $d_j$ be some points on $C$ 
between $z_j$ and the neighboring singularity in the clockwise and anticlockwise directions, respectively. 
In a neighborhood of $z_j$, let $g(s)$ be defined by the formula 
$\wh\phi_n(s^{-1})f(s)/(is)= |s-z_j|^{2\al_j}g(s)$.
We then obtain as in \ci{Khankel} for the part of (\ref{intdiff}) on the arc $(c_j,d_j)$:
\begin{multline}
\int_{c_j}^{d_j}{|s-z_j|^{2\al_j}\over s-z}g(s)ds-\lim_{\ep\to 0}\left[
 \left(\int_{c_j}^{z_j e^{-i\ep}}+\int_{z_j e^{i\ep}}^{d_j}\right)
{|s-z_j|^{2\al_j}\over s-z_j}g(s)ds-
{\ep^{2\al_j}\over 2\al_j}
(g(z_je^{-i\ep})-g(z_je^{i\ep}))\right]
\\
=\lim_{\ep\to 0}\frac{\pi (z-z_j)^{2\al_j}}{(zz_j)^{\al_j}\sin(2\pi\al_j)}
(e^{2\pi i\al_j-i\al\ell_R}g(z_je^{-i\ep})-
e^{-2\pi i\al_j-i\al\ell_L}g(z_je^{i\ep}))\\
+\al_j^{-1} O(z-z_j),\\
z\to z_j,\quad |z|>1,
\end{multline}
for $\al_j\neq 0, 1/2, 1,3/2,\dots$ (for $\al_j=1/2,1,3/2,\dots$ one obtains terms
involving $(z-z_j)^k\ln(z-z_j)$ vanishing as $z\to z_j$).
Here the constants $\ell_R$, $\ell_L$ depend on the 
choice of a branch for $(z-z_j)^{2\al}$ (whose cut is, recall, along the circle)
and their values will not be important below.  

Introduce a ``regularized'' version of the integral in a neighborhood of $z_j$:
\begin{multline}
{\int_{c_j}^{d_j}}^{(r)}{|s-z_j|^{2\al_j}\over s-z}g(s)ds\equiv
\int_{c_j}^{d_j}{|s-z_j|^{2\al_j}\over s-z}g(s)ds\\
-\lim_{\ep\to 0}\frac{\pi (z-z_j)^{2\al_j}}{(zz_j)^{\al_j}\sin(2\pi\al_j)}
(e^{2\pi i\al_j-i\al\ell_R}g(z_je^{-i\ep})-
e^{-2\pi i\al_j-i\al\ell_L}g(z_je^{i\ep})),
\end{multline}
for  $z$ in a complex neighborhood of  $z_j$ and  
$-1/2<\Re\al_j\le 0$, $\al_j\neq 0$. If $\Re\al_j>0$, we set the ``regularized'' integral
equal to the integral itself.

Denote by $\wt Y$ the matrix (\ref{RHM}), in which the integrals of the second
column are replaced by their ``regularized'' values in a neighborhood 
of each $z_j$.

Then, collecting our observations together, we can write (\ref{I1}) in the form
\be
I_1=
-\chi_n^{-1}{\partial\chi_n\over\partial \ga}
\left(n+\sum_{j=0}^m(\al_j+\bt_j)\right)+
\sum_{j=0}^m
\begin{cases}
{2\al_j z_j\over\chi_n}
{\partial\over\partial\ga}\left(\chi_n Y_{11}^{(n)}(z_j)\right)
\wt Y_{22}^{(n+1)}(z_j),& \al_j\neq 0\cr
{1\over 2\pi i}{\partial\phi_n(z_j)\over\partial \ga}\wh\phi_n(z_j^{-1})
\De f(z_j),& \al_j=0.
\end{cases}
\ee
where
\be\la{df}
\De f(z_j)=\lim_{\ep\to 0}(f(z_je^{-i\ep})-f(z_je^{i\ep})).
\ee

A similar analysis yields for $I_2$:
\be
I_2=
\chi_n^{-1}{\partial\chi_n\over\partial \ga}
\left(n+\sum_{j=0}^m(\al_j-\bt_j)\right)+
\sum_{j=0}^m
\begin{cases}
2\chi_n\al_j{\partial\over\partial\ga}\left(\chi_n^{-1}Y_{21}^{(n+1)}(z_j)\right)
\wt Y_{12}^{(n)}(z_j),& \al_j\neq 0\cr
{1\over 2\pi i}{\partial\wh\phi_n(z_j^{-1})\over\partial \ga}\phi_n(z_j)
\De f(z_j),& \al_j=0.
\end{cases}
\ee

Substituting these results into (\ref{mid2}) we obtain

\begin{proposition}\la{idDab}
Let $V(z)\equiv 0$.
Let $\ga=\al_k$ or $\ga=\bt_k$,  $k=0,1,\dots,m$,
and $D_n(f(z))\neq 0$ for all $n$. Then for any $n=1,2,\dots$, 
\begin{multline}\la{diffid1}
{\partial \over\partial \ga}\ln D_n(f(z))=
-2\chi_n^{-1}{\partial\chi_n\over\partial \ga}
\left(n+\sum_{j=0}^m\al_j\right)\\
+\sum_{j=0}^m
\begin{cases}
2\al_j\left\{
{\partial\over\partial\ga}\left(\chi_n Y_{11}^{(n)}(z_j)\right)
z_j\chi_n^{-1}\wt Y_{22}^{(n+1)}(z_j)-
{\partial\over\partial\ga}
\left(\chi_n^{-1}Y_{21}^{(n+1)}(z_j)\right)
\chi_n \wt Y_{12}^{(n)}(z_j)\right\},& \al_j\neq 0\cr
{1\over 2\pi i}\left\{
{\partial\phi_n(z_j)\over\partial \ga}\wh\phi_n(z_j^{-1})-
{\partial\wh\phi_n(z_j^{-1})\over\partial \ga}\phi_n(z_j)\right\}
\De f(z_j),& \al_j=0,
\end{cases}
\end{multline}
where $\De f(z_j)$ is defined in (\ref{df}).
\end{proposition}

%\begin{remark}
%The condition $D_n(f(z))\neq 0$ for all $n$ guarantees the 
%existence of the system of orthogonal polynomials corresponding to the weight
%$f(z)$ (see Section \ref{opuc}).
%\end{remark}

In Section \ref{Toep51} we substitute the asymptotics for $Y$ (found
in Section \ref{RHa}) in (\ref{diffid1})
and, by integrating, obtain part of Theorem \ref{asTop}
for $f(z)$ with $V(z)\equiv 0$, $|||\bt|||<1/2$. Further analysis of Section \ref{pFHext}
extends the result to  $|||\bt|||<1$. 

\begin{remark}\la{tau}
The differential identities (\ref{diffid1}) admit an interesting interpretation in the context
of the monodromy theory of the Fuchsian system of linear ODEs canonically related
to the Riemann-Hilbert problem (\ref{RHPYb})--(\ref{RHPYe}). We explain this connection 
in some detail in the Appendix. The results presented in the Appendix, however,
are not used in the main body of the paper.
\end{remark}
To extend the theorem to nonzero $V(z)$ we will use another differential 
identity. Let us introduce a parametric family of weights and the corresponding 
orthogonal  polynomials indexed by $t\in [0,1]$. Namely, let
\be\la{ft}
f(z,t)=(1-t+te^{V(z)})e^{-V(z)}f(z).
\ee
Thus $f(z,0)$ corresponds to $f(z)$ with $V=0$, whereas $f(z,1)$ gives the function 
(\ref{fFH}) we are interested in.

Note that
\be\la{ftdif}
{\partial f(z,t)\over\partial t}={ f(z,t)-f(z,0) \over t}.
\ee
Set now $\ga=t$ and replace the function $f(z)$ and the orthogonal 
polynomials in (\ref{mid}) by $f(z,t)$ and the polynomials orthogonal w.r.t. $f(z,t)$. 
Then the integral in the r.h.s. of (\ref{mid}) can be written as follows
(we assume $D_n\neq 0$ for all $n$):
\begin{multline}
{\partial\over \partial t}
\left[\frac{1}{2\pi}\int_0^{2\pi} 
\left(\phi_n(z,t){d\wh\phi_n(z^{-1},t)\over dz}
-\wh\phi_n(z^{-1},t){d\phi_n(z,t)\over dz}\right)z f(z,t) d\th
\right]\\
-\frac{1}{2\pi}\int_0^{2\pi} 
\left(\phi_n(z,t){d\wh\phi_n(z^{-1},t)\over dz}
-\wh\phi_n(z^{-1},t){d\phi_n(z,t)\over dz})\right)z
{ f(z,t)-f(z,0) \over t}d\th\\
={2n \over t}+
\frac{1}{2\pi t}\int_C 
\left(\phi_n(z,t){d\wh\phi_n(z^{-1},t)\over dz}
-\wh\phi_n(z^{-1},t){d\phi_n(z,t)\over dz}\right)z
f(z,0)dz.
\end{multline}
Therefore, we obtain
\begin{multline}\la{diffid2}
{\partial \over\partial t}\ln D_n(f(z,t))=
2n\left({1\over t}+\chi_n^{-1}{\partial\chi_n\over\partial t}\right)\\
+\frac{1}{2\pi t}\int_C 
\left(\phi_n(z,t){d\wh\phi_n(z^{-1},t)\over dz}
-\wh\phi_n(z^{-1},t){d\phi_n(z,t)\over dz}\right)z
f(z,0)dz.
\end{multline}
To write this identity in terms of the solution to the RHP (\ref{RHPYb}) - (\ref{RHPYc}),
note first that using the recurrence relation (see, e.g., Lemma 2.2 of \cite{DIKt1})
\[
\chi_n z^{-1}\widehat\phi_n(z^{-1})=
\chi_{n+1}\widehat\phi_{n+1}(z^{-1})-\widehat\phi_{n+1}(0)z^{-n-1}\phi_{n+1}(z),
\]
we have
\begin{equation}\label{Ytophi2}
Y_{21}(z,t) = - \chi_{n-1}z^{n-1}\wh{\phi}_{n-1}(z^{-1},t)
= - \chi_{n}z^{n}\wh{\phi}_{n}(z^{-1},t) +\wh{\phi}_{n}(0,t)
\phi_{n}(z,t).
\end{equation}
Now using the orthogonality relations 
(\ref{or0}) and the formulae (\ref{chidif1})
and  (\ref{chidif2}), we obtain from  (\ref{diffid2}) 

\begin{proposition}\la{idDt}
Let $f(z,t)$ be given by (\ref{ft}) and  $D_n(f(z,t))\neq 0$ for all $n$.
Let $\phi_k(z,t)$, $\wh\phi_k(z,t)$, 
$k=0,1,\dots$,
be the corresponding orthogonal polynomials.
Then for any $n=1,2,\dots$, 
\be\la{diffid2new}
{\partial \over \partial t}\ln D_n(f(z,t))=
\frac{1}{2\pi i}\int_{C}z^{-n} 
\left(Y_{11}(z,t){\partial Y_{21}(z,t)\over \partial z}
-Y_{21}(z,t){\partial Y_{11}(z,t)\over \partial z}\right)
{\partial f(z,t)\over \partial t}dz,
\ee
where the integration is over the unit circle.
\end{proposition}

%%%%%%%%%%%%%%%%% RHPa %%%%%%%%%%%%

\section{Asymptotics for the Riemann-Hilbert problem}\la{RHa}

The RHP of Section \ref{RHsection} was solved in \cite{DIKt1}.  
In this section we list the results from  \cite{DIKt1} we need below
for the proof of Theorem \ref{asTop}.
We always assume (for the rest of the paper) that $f(z)$ is given
by (\ref{fFH}) and that $\al_j\pm\bt_j\neq -1,-2,\dots$ for
all $j=0,1,\dots,m$.
In this section \ref{RHa} we also assume for simplicity that $z_0=1$
is a singularity. However, the results trivially extend to the case $\al_0=\bt_0=0$.
In this section \ref{RHa}, we further assume that $V(z)$ is analytic in a 
neighborhood of the unit circle.

First, set
\be\la{TY}
T(z)=Y(z)
\begin{cases}
z^{-n\si_3},& |z|>1\cr
I,& |z|<1.
\end{cases}
\ee

Now split the contour as shown in Figure 1. 
\begin{figure}
\centerline{\psfig{file=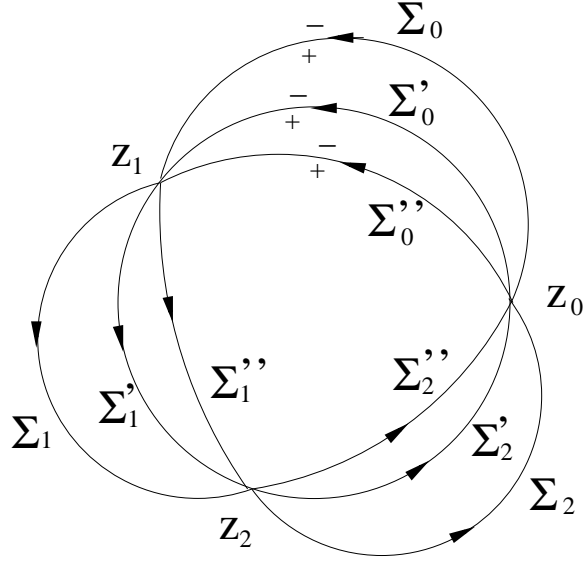,width=3.0in,angle=0}}
\vspace{0cm}
\caption{
Contour for the $S$-Riemann-Hilbert problem ($m=2$).}
\label{fig1}
\end{figure}
Set
\be\la{defS}
S(z)=
\begin{cases}T(z),& \mbox{for $z$ outside the lenses},\cr
T(z)\begin{pmatrix}1 & 0\cr f(z)^{-1}z^{-n}& 1\end{pmatrix}, & 
\mbox{for $|z|>1$ and inside the lenses},\cr
 T(z)\begin{pmatrix}1 & 0\cr -f(z)^{-1}z^n& 1\end{pmatrix},&
\mbox{for $|z|<1$ and inside the lenses}.
\end{cases}
\ee
Here $f(z)$ is the analytic continuation of $f(z)$ off the unit circle 
into the inside of the lenses as discussed following (\ref{abszzj}).

The function $S(z)$ satisfies the following Riemann Hilbert problem:
\begin{enumerate}
    \item[(a)]
        $S(z)$ is  analytic for $z\in\bbc \setminus\Sigma$, where 
$\Si=\cup_{j=0}^m(\Si_j\cup\Si'_j\cup\Si''_j)$.
    \item[(b)]  
The boundary values of $S(z)$ are related by the jump condition
\be\la{Sjump2}
   S_+(z) = S_-(z)
            \begin{pmatrix}
                  1 & 0\cr
                  f(z)^{-1}z^{\mp n } & 1\end{pmatrix},
            \qquad\mbox{$z\in\cup_{j=0}^{m}(\Si_j\cup\Si''_j)$},
\ee
where the minus sign in the exponent is on 
$\Si_j$, and plus on $\Si''_j$,
\be\la{Sjump1}
            S_+(z) = S_-(z)
              \begin{pmatrix}
                  0 & f(z)\cr
                  -f(z)^{-1} & 0\end{pmatrix},
            \qquad\mbox{$z\in\cup_{j=0}^{m}\Si'_j$.}
\ee
    \item[(c)]
$S(z)=I+O(1/z)$ as $z\to\infty$,
\item [(d)] 
As $z\to z_j$, $j=0,\dots,m$, $z\in\bbc\setminus C$ outside the lenses,
\be\la{Sd1}
S(z)=
\begin{pmatrix} O(1) & O(1)+O(|z- z_j|^{2\al_j})\cr 
O(1) & O(1)+O(|z- z_j|^{2\al_j})
\end{pmatrix}
\ee
if $\al_j\neq 0$, and
\be\la{Sd2}
S(z)=
\begin{pmatrix} O(1) & O(\ln|z- z_j|)\cr 
O(1) & O(\ln|z-z_j|)
\end{pmatrix}
\ee
if $\al_j=0$, $\bt_j\neq 0$.
The behavior of $S(z)$ for $z\to z_j$ in other sectors is obtained from these expressions by application of the appropriate jump conditions.
\end{enumerate}

We now present formulae for the parametrices which solve the model Riemann Hilbert problems
outside the neighborhoods $U_{z_j}$ of the points $z_j$, and inside those neighborhoods, respectively.
These parametrices match to the leading order in $n$ on the boundaries of the neighborhoods $U_{z_j}$,
and this matching allows us to construct the asymptotic solution to the RHP for $Y$.

The parametrix outside the $U_{z_j}$'s is the following:
\be\la{Ndef}
N(z)=
\begin{cases}
\mathcal{D}(z)^{\si_3},& |z|>1\cr
\mathcal{D}(z)^{\si_3}
\bm
0&1\cr -1&0
\em,
& |z|<1
\end{cases},\qquad z\in\bbc\setminus\cup_{j=0}^m U_{z_j}
\ee
where the Szeg\H o function
\be\la{48}
\mathcal{D}(z)=\exp {1\over 2\pi i}\int_C
{\ln f(s)\over s-z} ds,
\ee
is analytic away from the unit circle, and we have
\be\la{Dl1}
\mathcal{D}(z)
=e^{V_0}b_+(z)\prod_{k=0}^m\left({z-z_k\over z_k e^{i\pi}}\right)^{\al_k+\bt_k}
,\qquad
|z|<1.
\ee
and
\be\la{Dg1}
\mathcal{D}(z)
=b_-(z)^{-1}\prod_{k=0}^m\left({z-z_k\over z}\right)^{-\al_k+\bt_k}
,\qquad
|z|>1,
\ee
where $V_0$, $b_\pm(z)$ are defined in (\ref{WienH}).
Note that the branch of $(z-z_k)^{\pm\al_k+\bt_k}$ in (\ref{Dl1},\ref{Dg1}) 
is taken as discussed
following equation (\ref{abszzj}) above. In (\ref{Dg1}) for any $k$, the cut 
of the root $z^{-\al_k+\bt_k}$ is the line $\th=\th_k$ from 
$z=0$ to infinity, and $\th_k<\arg z<2\pi+\th_k$.

Inside each neighborhood $U_{z_j}$ the parametrix is given in terms of a
confluent hypergeometric function.
First, set
\be\la{zeta}
\zeta=n\ln{z\over z_j},
\ee
where $\ln x>0$ for $x>1$, and has a cut on the negative half of the real axis. 
Under this transformation
the neighborhood $U_{z_j}$ is mapped into a neighborhood of zero in the 
$\zeta$-plane. Note that the transformation $\zeta(z)$ is analytic, one-to-one,
and it takes an arc of the unit circle to an interval of the imaginary axis.
Let us now choose the exact form of the cuts $\Sigma$ in $U_{z_j}$ so that 
their images under the mapping $\zeta(z)$ are straight lines (Figure 2).

\begin{figure}
\centerline{\psfig{file=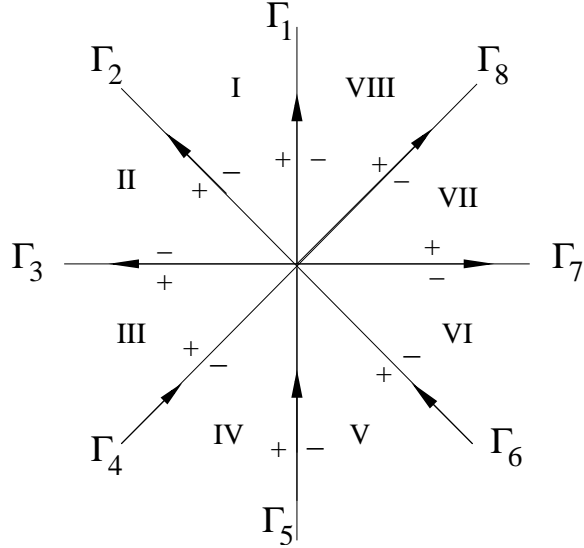,width=3.0in,angle=0}}
\vspace{0cm}
\caption{
The auxiliary contour for the parametrix at $z_j$.}
\label{fig2}
\end{figure}

We add one more jump contour to $\Sigma$ in $U_{z_j}$
which is the pre-image of the real line $\Ga_3$ and $\Ga_7$
in the $\zeta$-plane. 
This is needed below because of the non-analyticity of the function
$|z-z_j|^{\al_j}$.
Note that we can construct two different analytic continuations
of this function off the unit circle to the pre-images of the upper
and lower half $\ze$-plane, respectively.
Namely, let
\be
h_{\al_j}(z)=|z-z_j|^{\al_j},\qquad z=e^{i\th}
\ee
with the branches chosen as in (\ref{abszzj}).
As remarked above, (\ref{abszzj}) is continued analytically to neighborhoods of 
the arcs $0<\th<\th_j$, and
$\th_j<\th< 2\pi$. In $U_{z_j}$, we extend these neighborhoods to the pre-images 
of the lower and upper half $\ze$-plane (intersected with $\zeta(U_{z_j})$),
respectively. The cut of $h_{\al_j}$ is along 
the contours $\Ga_3$ and $\Ga_7$ in the $\ze$-plane.

For $z\to z_j$, $\zeta= n(z-z_j)/z_j+O((z-z_j)^2)$. 
We have $0<\arg\zeta<2\pi$, which follows from the choice of $\arg(z-z_j)$ in (\ref{abszzj}).

From now on we will provide the formulae for the parametrix only in the region $I$
(see \cite{DIKt1} for the complete results).
Set
\be\la{Fj}
F_j(z)
=e^{{V(z)\over 2}} \prod_{k=0}^m
\left(\frac{z}{z_k}\right)^{\bt_k/2} 
\prod_{k\neq j} h_{\al_k}(z) g_{\bt_k}(z)^{1/2}h_{\al_j}(z)
e^{-i\pi\al_j},\qquad \zeta\in I,
\quad z\in U_{z_j},\quad j\neq 0.
\ee
Note that this function is related to $f(z)$ as follows:
\be
F_j(z)^2 = f(z)e^{-2\pi i\al_j}g_{\bt_j}^{-1}(z)\qquad \ze\in I.
\ee
The functions $g_{\bt_k}(z)$ are defined in (\ref{g}). The formulae for 
$F_0(z)$ are the same, but with $g_{\bt_0}(z)$ replaced with
\be
\widehat g_{\bt_0}(z)=
\begin{cases}
e^{-i\pi\bt_0},& \arg z>0\cr
e^{i\pi\bt_0},& \arg z<2\pi
\end{cases},\qquad z\in U_{z_0}.
\ee

We then have the following expression for the parametrix $P_j(z)$ in the region 
$z(I)$ of $U_{z_j}$:
\be\la{Plb0}
P_{z_j}(z)=E(z)\Psi_j(\ze)F_j(z)^{-\si_3}z^{n\si_3/2},\qquad \ze\in I. 
\ee
Here
\be\la{E1}
E(z)=N(z)\ze^{\bt_j\si_3} F_j^{\si_3}(z) z_j^{-n\si_3/2}
\bm
e^{-i\pi(2\bt_j+\al_j)} &0\cr 0 & e^{i\pi(\bt_j+2\al_j)}
\em
\ee
and
\begin{multline}\label{PsiConfl1}
\Psi_j(\zeta)=\left(\begin{matrix}
\ze^{\al_j}\psi(\al_j+\bt_j,1+2\al_j,\ze)e^{i\pi(2\bt_j+\al_j)}e^{-\ze/2} \cr
-\ze^{-\al_j}
\psi(1-\al_j+\bt_j,1-2\al_j,\ze)e^{i\pi(\bt_j-3\al_j)}e^{-\ze/2}
{\Gamma(1+\al_j+\bt_j)\over\Gamma(\al_j-\bt_j)}
\end{matrix}
\right.\\
\left.
\begin{matrix}
-\ze^{\al_j}
\psi(1+\al_j-\bt_j,1+2\al_j,e^{-i\pi}\ze)e^{i\pi(\bt_j+\al_j)}e^{\ze/2}
{\Gamma(1+\al_j-\bt_j)\over\Gamma(\al_j+\bt_j)}
\cr
\ze^{-\al_j}
\psi(-\al_j-\bt_j,1-2\al_j,e^{-i\pi}\ze)e^{-i\pi\al_j}
e^{\ze/2}\end{matrix}\right),
\end{multline}
where $\psi(a,b,x)$ is the confluent hypergeometric function of the second kind, and 
$\Gamma(x)$ is Euler's $\Gamma$-function. 
Recall our assumption that $\al_j\pm\bt_j\neq -1,-2,\dots$

%Below we will need the behavior of $\psi(a,c,x)$ for small $x$. We have
%\begin{multline}\la{psi01}
%\psi(a,c,x)=
%\begin{cases}
%{\Ga(c-1)\over \Ga(a)}x^{1-c}\left(1+O(x\ln x)\right)+O(1),& \Re c>1\cr
%{\Ga(1-c)\over \Ga(1+a-c)}\left(1+O(x)\right)+
%{\Ga(c-1)\over \Ga(a)}x^{1-c},& \Re c=1, c\neq 1\cr
%-{1 \over \Gamma(a)}
%\left(\ln x +{\Gamma'(a)\over \Gamma(a)}-2C_\Gamma\right)+O(x\ln x),& c=1\cr
%{\Ga(1-c)\over \Ga(1+a-c)}\left(1+O(x\ln x)+O(x^{1-c})\right),& \Re c<1
%\end{cases},\qquad x\to 0,
%\end{multline}
%and
%\be\la{psi02}
%\psi(a,1,x)= -{1 \over \Gamma(a)}
%\left(\ln x +{\Gamma'(a)\over \Gamma(a)}-2C_\Gamma\right)+O(x\ln x),
%\qquad x\to 0,\qquad a \notin\{0,-1,-2,\dots\},
%\ee
%where $C_\Gamma=0.5772\dots$ is Euler's constant.

The matching condition for the parametrices $P_{z_j}$ and $N$ is the following for 
any $k=1,2,\dots$:
\be\la{Deas}
P_{z_j}(z)N^{-1}(z)=
I+\De_1(z)+\De_2(z)+\cdots+\De_{k}(z)+\De^{(r)}_{k+1},
\qquad z\in\partial U_{z_j}.
\ee
Every $\De_p(z)$, $\De^{(r)}_p(z)$, $p=1,2,\dots$, 
$z\in\partial U_{z_j}$ is of the form
\be\la{ordDe}
a_j^{-\si_3}O(n^{-p})a_j^{\si_3},\qquad a_j\equiv n^{\bt_j}
z_j^{-n/2}.
\ee
In particular, explicitly, on the part of $\partial U_{z_j}$ whose $\ze$-image is in $I$,
\be\la{De1}
\De_1(z)={1\over\ze}
\bm
-(\al_j^2-\bt_j^2)&
{\Ga(1+\al_j+\bt_j)\over \Ga(\al_j-\bt_j)}
\left({\mathcal{D}(z)\over\ze^{\bt_j}F_j(z)}\right)^2
z_j^n e^{i\pi(2\bt_j-\al_j)}\cr
-{\Ga(1+\al_j-\bt_j)\over \Ga(\al_j+\bt_j)}
\left({\mathcal{D}(z)\over\ze^{\bt_j}F_j(z)}\right)^{-2}
z_j^{-n} e^{-i\pi(2\bt_j-\al_j)}&
\al_j^2-\bt_j^2
\em,
\ee
which extends to a meromorphic function in 
a neighborhood of $U_{z_j}$ with a simple pole at $z=z_j$.

The error term $\De^{(r)}_{k+1}$ in (\ref{Deas}) is uniform in $z$ on $\partial U_{z_j}$.

At the point $z_j$ we have
\be\la{etab}
F_j(z)=\eta_j e^{-3i\pi\al_j/2}z_j^{-\al_j}
u^{\al_j}(1+O(u)),\qquad u=z-z_j,\qquad \ze\in I,
\ee
where
\be\la{eta}
\eta_j=e^{V(z_j)/2}\exp\left\{-{i\pi\over 2}\left(\sum_{k=0}^{j-1}\bt_k-
\sum_{k=j+1}^m\bt_k\right)\right\}
\prod_{k\neq j}\left({z_j\over z_k}\right)^{\bt_k/2}|z_j-z_k|^{\al_k},
\ee
and
\begin{eqnarray}\la{mub}
\left({\mathcal{D}(z)\over \ze^{\bt_j}F_j(z)}\right)^2=\mu^2_j e^{i\pi(\al_j-2\bt_j)}
n^{-2\bt_j}(1+O(u)),
\qquad u=z-z_j,\qquad \ze\in I,
\\
\mu_j=\left(e^{V_0}\frac{b_+(z_j)}{b_-(z_j)}\right)^{1/2}
\exp\left\{-{i\pi\over2}\left(\sum_{k=0}^{j-1}\al_k-
\sum_{k=j+1}^m\al_k\right)\right\}
\prod_{k\neq j}\left({z_j\over z_k}\right)^{\al_k/2}|z_j-z_k|^{\bt_k}.
\la{mu}
\end{eqnarray}
The sums from $0$ to $-1$ for $j=0$ and from $m+1$ to $m$ for $j=m$ are set to zero.

%%%%%%%%%%% R %%%%%%%%%%%%%

\subsection{R-RHP}\la{RRHP}
Let
\be
R(z)=\begin{cases}S(z)N^{-1}(z),& 
z\in U_\infty\setminus\Gamma,\qquad U_\infty=\bbc\setminus\cup_{j=0}^m U_{z_j},\cr
S(z)P_{z_j}^{-1}(z),& 
z\in U_{z_j}\setminus\Gamma,\qquad j=0,\dots, m.
\end{cases}\la{wtR}
\ee
It is easy to verify that this function has jumps only on 
$\partial U_{z_j}$, and the parts of 
$\Si_j$,  $\Si^{''}_j$
lying outside the neighborhoods $U_{z_j}$ (we
denote these parts without the end-points $\Si^\mathrm{out}$, $\Si^{''\mathrm{out}}$). 
The full contour $\Gamma$ is shown in Figure 3. Away from $\Gamma$, as a standard argument
shows, $R(z)$ is analytic.
Moreover, we have: $R(z)=I+O(1/z)$ as $z\to\infty$. 

\begin{figure}
\centerline{\psfig{file=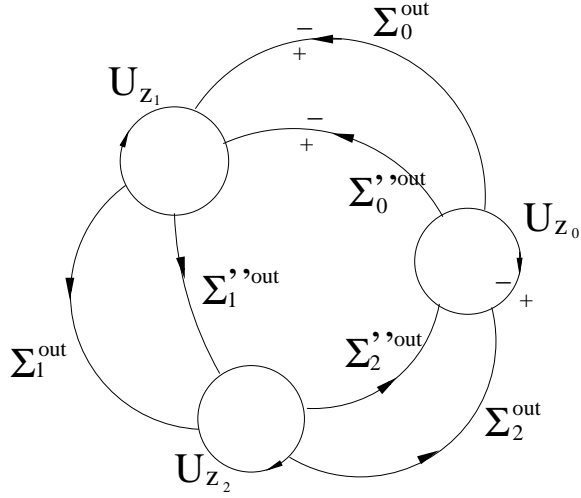,width=3.0in,angle=0}}
\vspace{0cm}
\caption{ 
Contour $\Gamma$ for the $R$ and $\wt R$ Riemann-Hilbert problems ($m=2$).}
\label{fig3}
\end{figure}

The jumps of $R(z)$ are as follows:
\begin{eqnarray}
R_+(z)&=&R_-(z)N(z)
\bm
1&0\cr f(z)^{-1}z^{-n}&1
\em
N(z)^{-1},\qquad z\in\Si_j^\mathrm{out},\la{Rs1}\\
R_+(z)&=&R_-(z)N(z)
\bm
1&0\cr f(z)^{-1}z^{n}&1
\em
N(z)^{-1},\qquad z\in\Si_j^{''\mathrm{out}},\la{Rs2}\\
R_+(z)&=&R_-(z)P_{z_j}(z)N (z)^{-1},\qquad 
z\in \partial U_{z_j}\setminus
\mbox{\{ intersection points\}},\la{Rs3}\\
j=0,\dots,m.\nonumber
\end{eqnarray}

The jump matrix on $\Si^\mathrm{out}$, $\Si^{''\mathrm{out}}$ can be 
estimated uniformly in $\al_j$, $\bt_j$ as $I+O(\exp(-\ep n))$, where 
$\ep$ is a positive constant.
The jump matrices on $\partial U_{z_j}$ admit a uniform expansion (\ref{Deas}) in 
inverse powers of $n$ conjugated by $n^{\bt_j\si_3}z_j^{-n\si_3/2}$,
and (\ref{ordDe}) is of order $n^{2\max_j|\Re\bt_j|-p}$.
To obtain the standard solution of the $R$-RHP in terms of a Neuman series
(see, e.g., \cite{Dstrong}) we must have 
$n^{2\max_j|\Re\bt_j|-1}=o(1)$, that is
$\Re\bt_j\in (-1/2,1/2)$ for all $j=0,1,\dots,m$.
However, it is possible to obtain the solution in any half-closed or open interval of length 1,
i.e. for $|||\bt|||<1$, as follows.
 
Let $|||\bt|||<1$ and consider the transformation 
\be\la{Rtilde}
\wt R(z)=n^{\om\si_3}R(z) n^{-\om\si_3}\qquad z\in\bbc\setminus\Gamma,
\ee
where
\be\la{omega}
\om={1\over 2}(\min_j\Re\bt_j+\max_j\Re\bt_j)
\ee
which ``shifts'' all $\Re\bt_j$ (in the conjugation $n^{\bt_j}$ terms of (\ref{ordDe}) for
the jump matrix in (\ref{Rs3})) into the interval $(-1/2,1/2)$.  
Note that $\om=\Re\bt_{j_0}$ if only one $\Re\bt_{j_0}\neq 0$, and $\om=0$
if all $\Re\bt_j=0$.

Now in the RHP for $\wt R(z)$, the condition at infinity and the uniform 
exponential estimate $I+O(\exp(-\ep n))$ (with different $\ep$)
of the jump matrices on $\Si^\mathrm{out}$, $\Si^{''\mathrm{out}}$ is preserved,
while the jump matrices on $\partial U_{z_j}$ have the form:
\be\la{Deas2}
I+n^{\om\si_3}\De_1(z)n^{-\om\si_3}+\cdots+
n^{\om\si_3}\De_{k}(z)n^{-\om\si_3}+
n^{\om\si_3}\De_{k+1}^{(r)}(z)n^{-\om\si_3},\qquad 
z\in\partial U_{z_j},
\ee
where the order of each $n^{\om\si_3}\De_p(z)n^{-\om\si_3}$, 
$n^{\om\si_3}\De_p^{(r)}(z)n^{-\om\si_3}$, $p=1,2,\dots$, 
$z\in\cup_{j=0}^m\partial U_{z_j}$
is 
\[
O(n^{2\max_j|\Re\bt_j-\om|-p})=O(n^{|||\bt|||-p}).
\]
This implies that the standard analysis can be applied to the $\wt R$-RHP problem
in the range $\Re\bt_j\in (q-1/2,q+1/2)$, $j=0,1,\dots,m$,
for any $q\in\bbr$, and we obtain the asymptotic expansion
\be\la{Ras}
\wt R(z)=I+\sum_{p=1}^{k}\wt R_p(z)+\wt R^{(r)}_{k+1}(z),\quad p=1,2\dots
\ee
uniformly for all $z$ and for $\bt_j$ in  bounded sets 
of the strip $q-1/2<\Re\bt_j< q+1/2$, $j=0,1,\dots m$, i.e. $|||\bt|||<1$, 
provided
$\al_j\pm\bt_j$ are outside neighborhoods of the points $\al_j\pm\bt_j=-1,-2,\dots$
(cf. (\ref{De1})).
 
The functions $\wt R_j(z)$ are computed recursively. 
We will need explicit expressions only for 
the first two.  The first one is found from the conditions
that $\wt R_1(z)$ is analytic outside 
$\partial U=\cup_{j=0}^m\partial U_{z_j}$, 
$\wt R_1(z)\to 0$ as $z\to \infty$, and 
\be
\wt R_{1,+}(z)=\wt R_{1,-}(z)+n^{\om\si_3}\De_1(z)n^{-\om\si_3},\qquad 
z\in\partial U.
\ee
The solution is easily found. First denote
\be
R_k(z)\equiv n^{-\om\si_3}\wt R_k(z)n^{\om\si_3},\qquad
R_k^{(r)}(z)\equiv n^{-\om\si_3}\wt R_k^{(r)}(z)n^{\om\si_3},\qquad
k=1,2,\dots,
\ee
and write for $R$:
\begin{multline}
R_1(z)=
{1\over 2\pi i}\int_{\partial U}
{\De_1(x)dx\over x-z}\\
=\begin{cases}
\sum_{k=0}^m {A_k\over z-z_k},& z\in\bbc\setminus \cup_{j=0}^m U_{z_j}\cr
\sum_{k=0}^m {A_k\over z-z_k}-\De_1(z),& z\in U_{z_j},\quad j=0,1,\dots, m.
\end{cases},\qquad 
\partial U=\cup_{j=0}^m\partial U_{z_j}.
\la{plem}
\end{multline}
where the contours in the integral are traversed in the negative direction,
and $A_k$ are the coefficients in the Laurent expansion of $\De_1(z)$:
\be
\De_1(z)={A_k\over z-z_k}+B_k+O(z-z_k),\qquad z\to z_k,\qquad k=0,1,\dots,m.
\ee
The coefficients are easy to compute using (\ref{Deas}) and (\ref{mu}):
\be\la{A}
A_k\equiv A^{(n)}_k={z_k\over n}
\bm
-(\al_k^2-\bt_k^2)&
{\Ga(1+\al_k+\bt_k)\over \Ga(\al_k-\bt_k)}
z_k^n \mu^2_k n^{-2\bt_k}\cr
-{\Ga(1+\al_k-\bt_k)\over \Ga(\al_k+\bt_k)}
z_k^{-n} \mu_k^{-2}n^{2\bt_k}&
\al_k^2-\bt_k^2
\em.
\ee
The function $\wt R_2$ is now found from the conditions that $\wt R_2(z)\to 0$ 
as $z\to\infty$, is analytic outside $\partial U$, and
\be
\wt R_{2,+}(z)=\wt R_{2,-}(z)+\wt R_{1,-}(z)n^{\om\si_3}\De_1(z)n^{-\om\si_3}+
n^{\om\si_3}\De_2(z)n^{-\om\si_3} ,\qquad 
z\in\partial U.
\ee
The solution to this RHP is
\be\la{R2}
\wt R_2(z)=
 {1\over 2\pi i}\int_{\partial U}\left(
\wt R_{1,-}(x)n^{\om\si_3}\De_1(x)n^{-\om\si_3}+
n^{\om\si_3}\De_2(x)n^{-\om\si_3}
\right){dx\over x-z}.
\ee
At the $k$'th step we have the RHP for $\wt R_k(z)$ with the same analyticity 
condition and the condition at infinity, and the following jump:
\be\la{jumpRk}
\wt R_{k,+}(z)=\wt R_{k,-}(z)+
\sum_{p=1}^k
\wt R_{k-p,-}(z)n^{\om\si_3}\De_{p}(z)n^{-\om\si_3} ,\qquad 
z\in\partial U,
\ee
where $\wt R_0(z)\equiv I$.

We will now discuss the way in which the general $\wt R_k(z)$ depends
on $n$. In particular, we will discuss its order in $n$.
First note that
\begin{eqnarray}
\wt R_1(z)\sim {1\over n}
\bm
1& \sum_j b_j^{-2}\cr
\sum_j b_j^{2} & 1
\em,\qquad
\wt R_2(z)\sim {1\over n^2}
\bm
1+\de' n^2& \sum_j b_j^{-2}\cr
\sum_j b_j^{2} & 1+\de' n^2
\em,\\
b_j\equiv n^{\bt_j-\om}z_j^{-n/2},\qquad 
\de'\sim \sum_{j,k}n^{2((\bt_j-\bt_k)-1)}\left({z_k\over z_j}\right)^n.
\end{eqnarray}
Here the notation $A\sim B$ means $A'=B'$, where $X'$ is $X$ in which each matrix 
element and each term in the sums is multiplied by a suitable constant 
{\it independent} of $n$.
Starting with these expressions, and noting from (\ref{jumpRk}) that
\be
\wt R_k(z)\sim \sum_{p=1}^k
\wt R_{k-p,-}(z)n^{\om\si_3}\De_{p}(z)n^{-\om\si_3}\sim
\wt R_{k-1,-}(z)n^{\om\si_3}\De_{1}(z)n^{-\om\si_3}+{1\over n}\wt R_{k-1,-}(z),
\ee
we obtain by induction:
\begin{eqnarray}
&\wt R_{2p+1}(z)\sim
{1\over n^{2p+1}}\sum_{k=0}^p \left(\de' n^2\right)^k
\bm
1& \sum_j b_j^{-2}\cr
\sum_j b_j^{2} & 1
\em,\\
&\wt R_{2p+2}(z)\sim
{1\over n^{2p+2}}\sum_{k=0}^p \left(\de' n^2\right)^k
\bm
1+\de' n^2 & \sum_j b_j^{-2}\cr
\sum_j b_j^{2} & 1+\de' n^2
\em,\qquad p=0,1,\dots.
\end{eqnarray}
In particular,
\begin{eqnarray}\la{ordR}
&\wt R_{2p+1}(z)={{\de'}^p\over n}
O\bm
1& \sum_j b_j^{-2}\cr
\sum_j b_j^{2} & 1
\em,\qquad
\wt R_{2p+2}(z)={{\de'}^p\over n^2}
O\bm
1+\de' n^2 & \sum_j b_j^{-2}\cr
\sum_j b_j^{2} & 1+\de' n^2
\em,\\
&O(\de')=O(\de),\quad \de=\max_{j,k} n^{2(\Re(\bt_j-\bt_k)-1)}=n^{2(|||\bt|||-1)},
\qquad p=0,1,\dots,
\end{eqnarray}
as $n\to\infty$.
Here $O(A)$  represent $2\times 2$ matrices with elements
of the corresponding order.

Finally, note that the error term in (\ref{Ras}) is
\be\la{orderr}
\wt R^{(r)}_k(z)=O(|\wt R_{k}(z)|+|\wt R_{k+1}(z)|).
\ee
  
In particular, as is clear from the above, if there is only one 
nonzero $\bt_{j_0}$, we obtain the expansion 
purely in inverse integer powers of $n$ valid in fact for all $\bt_{j_0}\in \bbc$
uniformly in bounded sets of the complex plane.

It is clear from the construction and the properties of the asymptotic series
of the confluent hypergeometric function that
the error terms $\wt R^{(r)}_k(z)$ are uniform for $\bt_j$ in bounded subsets of the strip
$q-1/2<\Re\bt_j< q+1/2$, $j=0,1,\dots m$, for $\al_j$ in bounded sets of the half-plane
$\Re\al_j>-1/2$, and for $\al_j\pm\bt_j$ away from neighborhoods of the negative integers.
Moreover, the series (\ref{Ras}) is differentiable in $\al_j$, $\bt_j$.

\section{Asymptotics for differential identities and integration. 
Proof of Theorem \ref{asTop}}\la{Toep}
\subsection{Pure Fisher-Hartwig singularities. The case $|||\bt|||<1/2$.}\la{Toep51}
First, we will prove the theorem for $V(z)\equiv 0$ and 
$|||\bt|||=\max_{j,k}|\Re\bt_j-\Re\bt_k|<1/2$.
The proof is based on the differential identity (\ref{diffid1}). First, we show that (\ref{diffid1})
has the following asymptotic form.

\begin{proposition}\la{asid} Let $(\al_0,\bt_0,\dots,\al_m,\bt_m)$ be in a compact subset,
denote it $\Lambda$,
belonging to the subset $|||\bt|||<1$, $\al_j\pm\bt_j\neq -1,-2,\dots$ of the parameter space 
${\mathcal P}=\{(\al_0,\bt_0,\dots,\al_m,\bt_m):\, \al_j,\bt_j\in\bbc,\, \Re\al_j>-1/2\}$
and including the point $\al_j=\bt_j=0$, $j=0,1,\dots,m$.
Let $\bt_j=0$ if $\al_j=0$,  $j=0,1,\dots,m$, $\de=n^{2(|||\bt|||-1)}$. 
Then for $n\to\infty$, and $\nu=0,1,\dots, m$,
\begin{multline}\la{asidal}
{\partial\over\partial\al_\nu}\ln D_n(f(z))=
2\al_\nu\\
+(\al_\nu+\bt_\nu)\left[
{\partial\over\partial\al_\nu}\ln
{\Ga(1+\al_\nu+\bt_\nu)\over\Ga(1+2\al_\nu)}+\ln n
\right]+
(\al_\nu-\bt_\nu)\left[
{\partial\over\partial\al_\nu}\ln
{\Ga(1+\al_\nu-\bt_\nu)\over\Ga(1+2\al_\nu)}+\ln n
\right]\\
-\sum_{j\neq\nu}\left[
(\al_j+\bt_j)\ln{z_j-z_\nu \over z_j}+
(\al_j-\bt_j)\ln{z_j-z_\nu \over z_\nu e^{i\pi}}\right]
+2\pi i\sum_{j=0}^{\nu-1}(\al_j+\bt_j)+O(n^{-1}\ln n)+O(\de n\ln n)
\end{multline}
and
\begin{multline}\la{asidbt}
{\partial\over\partial\bt_\nu}\ln D_n(f(z))=
-2\bt_\nu\\
+(\al_\nu+\bt_\nu)\left[
{\partial\over\partial\bt_\nu}\ln\Ga(1+\al_\nu+\bt_\nu)-\ln n
\right]+
(\al_\nu-\bt_\nu)\left[
{\partial\over\partial\bt_\nu}\ln\Ga(1+\al_\nu-\bt_\nu)+\ln n
\right]\\
+\sum_{j\neq\nu}\left[
(\al_j+\bt_j)\ln{z_j-z_\nu \over z_j}-
(\al_j-\bt_j)\ln{z_j-z_\nu \over z_\nu e^{i\pi}}\right]
-2\pi i\sum_{j=0}^{\nu-1}(\al_j+\bt_j)+O(n^{-1}\ln n)+O(\de n\ln n).
\end{multline}
\end{proposition}

\begin{remark}
The error term $O(\de n\ln n)=o(1)$ uniformly in $\Lambda$ if $|||\bt|||\le 1/2-\ep$, $\ep>0$.
In fact, the estimate for the error term can be considerably improved: see next section.
\end{remark}

\begin{remark}
The case $\bt_j=0$ if $\al_j=0$ is all we need below. After the proof of Proposition \ref{asid},
we will integrate the identity (\ref{asidal}) to obtain the asymptotics of $D_n$ for all
$\al_j\neq 0$, $\bt_j=0$. We then integrate the identity (\ref{asidbt}) and obtain the asymptotics
of $D_n$ for all $\al_j\neq 0$, $\bt_j\neq 0$. This gives the general result for $V\equiv 0$, $|||\bt|||<1/2$,
since we can set any $\al_j=0$ using the uniformity of the asymptotic expansion in the $\al_j$'s.
\end{remark}

\begin{proof}
Assume that for all $j$, $\bt_j=0$ if $\al_j=0$, and $D_k(f)\neq 0$, $k=1,2,\dots$.
Then we can rewrite (\ref{diffid1}) in the form:
\begin{multline}\la{diffid3}
{\partial\over\partial\ga}\ln D_n(f(z))=
-2{{\partial\chi_n\over\partial\ga}\over\chi_n}
\left(n+\sum_{j=0}^m\al_j\left\{1-\wt Y_{12}^{(n)}(z_j)Y_{21}^{(n+1)}(z_j)-
z_jY_{11}^{(n)}(z_j)\wt Y_{22}^{(n+1)}(z_j)\right\}\right)\\
-
2\sum_{j=0}^m\al_j\left(\wt Y_{12}^{(n)}(z_j){\partial\over\partial\ga}Y_{21}^{(n+1)}(z_j)-
z_j{\partial\over\partial\ga}Y_{11}^{(n)}(z_j)\wt Y_{22}^{(n+1)}(z_j)
\right),
\end{multline}
where $\ga=\al_j$ or $\ga=\bt_j$.
We now estimate the
right-hand side of this identity as $n\to\infty$. 
The asymptotics of $\chi_n$ were found in (\ci{DIKt1}, Theorem 1.8). We need these asymptotics
here in the case $V\equiv 0$:
\begin{multline}\la{aschi}
\chi_{n-1}^2=
1-{1\over n}\sum_{k=0}^m (\al_k^2-\bt_k^2)\\
+\sum_{j=0}^m\sum_{k\neq j}{z_k\over z_j-z_k}
\left({z_j\over z_k}\right)^n n^{2(\bt_k-\bt_j-1)}{\nu_j\over\nu_k}
{\Ga(1+\al_j+\bt_j)\Ga(1+\al_k-\bt_k)\over
\Ga(\al_j-\bt_j)\Ga(\al_k+\bt_k)}\\
+ O(\de^2)+O(\de/n),\qquad  \de=n^{2(|||\bt|||-1)}, \qquad n\to\infty,
\end{multline}
where
\be\la{nu}
\nu_j=\exp\left\{-i\pi\left(\sum_{p=0}^{j-1}\al_p-
\sum_{p=j+1}^m\al_p\right)\right\}
\prod_{p\neq j}\left({z_j\over z_p}\right)^{\al_p}
|z_j-z_p|^{2\bt_p}.
\ee

The asymptotics of $\wt Y^{(n)}(z_j)$ were also found in \ci{DIKt1} (Eq. (7.11)-(7.21)). Namely,  
\be\la{Yend}
\wt Y^{(n)}(z_j)=(I+r^{(n)}_j)L^{(n)}_j,\qquad
L^{(n)}_j=\bm
M_{21}\mu_j\eta_j^{-1} n^{\al_j-\bt_j} z_j^{n} &
M_{22}\mu_j\eta_j n^{-\al_j-\bt_j}\cr
-M_{11}\mu^{-1}_j \eta_j^{-1} n^{\al_j+\bt_j} & 
-M_{12}\mu^{-1}_j \eta_j n^{-\al_j+\bt_j} z_j^{-n}
\em,
\ee
where $r_j=R_1^{(r)}(z_j)$, the parameters $\eta_j$, $\mu_j$ are given by (\ref{eta}), (\ref{mu}), and
\[
M=
\bm
{\Ga(1+\al_j-\bt_j)\over \Ga(1+2\al_j)}&
 -{\Ga(2\al_j)\over \Ga(\al_j+\bt_j)}\cr
{\Ga(1+\al_j+\bt_j)\over \Ga(1+2\al_j)}&
 {\Ga(2\al_j)\over \Ga(\al_j-\bt_j)}
\em.
\]

Note that the matrix $L^{(n)}_j$ has the structure
\be\la{LhatL}
L^{(n)}_j=n^{-\bt_j\si_3}\wh L^{(n)}_j n^{\al_j\si_3},
\ee
where $\wh L$ depends on $n$ only via the oscillatory terms $z_j^n$.

Let us now obtain the asymptotics for the following combination appearing in
(\ref{diffid3}):
\begin{multline}
{\partial\over\partial\ga}Y_{11}^{(n)}(z_j)\wt Y_{22}^{(n+1)}(z_j)
=
{\partial\over\partial\ga}\left((1+r^{(n)}_{11})L^{(n)}_{11}+
r^{(n)}_{12}L^{(n)}_{21}\right)
\left(r^{(n+1)}_{21}L^{(n+1)}_{12}+(1+r^{(n+1)}_{22})L^{(n+1)}_{22}
\right)\\
=
\left(L^{(n)}_{11}L^{(n+1)}_{22}(1+r^{(n+1)}_{22})+
L^{(n)}_{11}L^{(n+1)}_{12}r^{(n+1)}_{21}\right)
\left[(1+r^{(n)}_{11}){\partial\over\partial\ga}\ln L^{(n)}_{11}+
{\partial\over\partial\ga}r^{(n)}_{11}\right]\\
+
\left(L^{(n)}_{21}L^{(n+1)}_{22}(1+r^{(n+1)}_{22})+
L^{(n)}_{21}L^{(n+1)}_{12}r^{(n+1)}_{21}\right)
\left[r^{(n)}_{12}{\partial\over\partial\ga}\ln L^{(n)}_{21}+
{\partial\over\partial\ga}r^{(n)}_{12}\right].
\end{multline}
We omit the lower index $j$ of $r$ and $L$ for simplicity of notation.

Now the explicit formula for $L$ and the estimates for $\wt R(z)$ imply that
\begin{eqnarray}
&L^{(n)}_{11}L^{(n+1)}_{22}=O(1),\qquad L^{(n)}_{21}L^{(n+1)}_{12}=O(1),\\
&L^{(n)}_{11}L^{(n+1)}_{12}r^{(n+1)}_{21}=O(n\de),\qquad
L^{(n)}_{21}L^{(n+1)}_{22}r^{(n)}_{12}=O(n\de),\qquad
r^{(n+1)}_{21}r^{(n)}_{12}=O(\de),\\
&{\partial\over\partial\ga}\ln L^{(n)}=O(\ln n),\qquad
{\partial\over\partial\ga}r=O(r)\ln n,
\end{eqnarray}
where as before
\[
\de=n^{2(|||\bt|||-1)}.
\]
Therefore,
\begin{multline}\la{1}
{\partial\over\partial\ga}Y_{11}^{(n)}(z_j)\wt Y_{22}^{(n+1)}(z_j)
=
L^{(n)}_{11}L^{(n+1)}_{22}{\partial\over\partial\ga}\ln L^{(n)}_{11}+
O\left({\ln n\over n}\right)+O(\de n\ln n)\\
={\al_j+\bt_j \over 2\al_j z_j}{\partial\over\partial\ga}\ln L^{(n)}_{11}+
O\left({\ln n\over n}\right)+O(\de n\ln n).
\end{multline}
Similarly, we obtain
\be\la{2}
\wt Y_{12}^{(n)}(z_j){\partial\over\partial\ga}Y_{21}^{(n+1)}(z_j)=
-{\al_j-\bt_j \over 2\al_j}{\partial\over\partial\ga}\ln L^{(n+1)}_{21}+
O\left({\ln n\over n}\right)+O(\de n\ln n),
\ee
and furthermore,
\be\la{3}
\wt Y_{12}^{(n)}(z_j)Y_{21}^{(n+1)}(z_j)=O(\de n)+O(1),\qquad
Y_{11}^{(n)}(z_j)\wt Y_{22}^{(n+1)}(z_j)=O(\de n)+O(1).
\ee
Note that because of the special structure of (\ref{LhatL}),
the quantity $n^{\al_j}$ does not appear in any of the products
(\ref{1})--(\ref{3}). Substituting (\ref{1})--(\ref{3}) into (\ref{diffid3})
and using the asymptotics (\ref{aschi}), we obtain
\begin{multline}\la{diffidmid}
{\partial\over\partial\ga}\ln D_n(f(z))=
{\partial\over\partial\ga}\left[\sum_{j=0}^m(\al_j^2-\bt_j^2)\right]\\
+\sum_{j=0}^m\left[ (\al_j+\bt_j){\partial\over\partial\ga}\ln L^{(n)}_{j,11}+
(\al_j-\bt_j){\partial\over\partial\ga}\ln L^{(n+1)}_{j,21}
\right]+O\left({\ln n\over n}\right)+O(\de n\ln n).
\end{multline}

Let us calculate the logarithmic derivatives appearing in (\ref{diffidmid}). From (\ref{Yend}), (\ref{mu}),
and (\ref{eta}) it is easy to obtain for the derivatives w.r.t. $\al_\nu$, 
$\nu=0,1,\dots,m$:
\begin{eqnarray}
&{\partial\over\partial\al_\nu}\ln L^{(n)}_{\nu,11}=
{\partial\over\partial\al_\nu}\ln
{\Ga(1+\al_\nu+\bt_\nu)\over \Ga(1+2\al_j)}+\ln n;\\
&{\partial\over\partial\al_\nu}\ln L^{(n)}_{j,11}=
-\ln{z_j-z_\nu \over z_j}+2\pi i,\quad j<\nu;\qquad
{\partial\over\partial\al_\nu}\ln L^{(n)}_{j,11}=
-\ln{z_j-z_\nu \over z_j},\quad j>\nu;\\
&{\partial\over\partial\al_\nu}\ln L^{(n+1)}_{\nu,21}=
{\partial\over\partial\al_\nu}\ln
{\Ga(1+\al_\nu-\bt_\nu)\over \Ga(1+2\al_j)}+\ln n;\\
&{\partial\over\partial\al_\nu}\ln L^{(n+1)}_{j,21}=
-\ln{z_j-z_\nu \over z_\nu e^{i\pi}},\quad j\neq\nu.
\end{eqnarray}
Similarly, we obtain for the derivatives w.r.t. $\bt_\nu$:
\begin{eqnarray}
&{\partial\over\partial\bt_\nu}\ln L^{(n)}_{\nu,11}=
{\partial\over\partial\bt_\nu}\ln\Ga(1+\al_\nu+\bt_\nu)-\ln n;\\
&{\partial\over\partial\bt_\nu}\ln L^{(n)}_{j,11}=
\ln{z_j-z_\nu \over z_j}-2\pi i,\quad j<\nu;\qquad
{\partial\over\partial\bt_\nu}\ln L^{(n)}_{j,11}=
\ln{z_j-z_\nu \over z_j},\quad j>\nu;\\
&{\partial\over\partial\bt_\nu}\ln L^{(n+1)}_{\nu,21}=
{\partial\over\partial\bt_\nu}\ln
\Ga(1+\al_\nu-\bt_\nu)+\ln n;\\
&{\partial\over\partial\bt_\nu}\ln L^{(n+1)}_{j,21}=
-\ln{z_j-z_\nu \over z_\nu e^{i\pi}},\quad j\neq\nu.\la{mark}
\end{eqnarray}
Combining these results with (\ref{diffidmid}) we obtain (\ref{asidal}), (\ref{asidbt})
on condition that $D_k(f)\neq 0$, $k=1,2,\dots$, or equivalently (see Section \ref{RHsection}),
$(\al_0,\bt_0,\dots,\al_m,\bt_m)\in\Lambda\setminus\Om_{k_0}$. 
This condition can be replaced simply by $(\al_0,\bt_0,\dots,\al_m,\bt_m)\in\Lambda$ in the following way. 
Let
\be\la{cond0}
\bt_0=\al_1=\cdots=\al_m=\bt_m=0.
\ee
Let $\Om_{k_0}(\al_0)$ be the subset of $\Om_{k_0}$ with $\al_j$, $\bt_j$ fixed by (\ref{cond0}).
Since $D_k(f)\equiv D_k(f(\al_0;z))$ is an analytic function of $\al_0$ and $D_k(1)\neq 0$, the set  
$\Om_{k_0}(\al_0)$ is finite. 
Let us rewrite the identity (\ref{asidal}) with $\nu=0$ and assuming (\ref{cond0}) in 
the form $H'(\al_0)=0$, where $H(\al_0)=D_n(f(\al_0;z))\exp(-\int_0^{\al_0} r(n,s)ds)$ and
where $r(n,\al_0)$ is the r.h.s. of (\ref{diffid1}) with $\ga=\al_0$ and assuming (\ref{cond0}).
Since the expression (\ref{asidal}) for $r(n,\al_0)$ holds uniformly and is continuous 
for $\al_0\in\Lambda$ provided $n$ is larger than some $k_0(\Lambda)$,
and $D_n(f(\al_0;z))$ and its derivative are continuous, the function
$H(\al_0)$ is continuously differentiable for all $n>k_0(\Lambda)$. Hence,
$H'(\al_0)=0$ for {\it all} $\al_0\in\Lambda$ and $n>k_0(\Lambda)$. 
Taking into account that $H(0)=D_n(1)\neq 0$,
we conclude that  $D_n(f(\al_0;z))$ is nonzero, and that the identity (\ref{asidal}) under (\ref{cond0})
is, in fact, true for all $\al_0\in\Lambda$ if $n$ is sufficiently large
(larger than $k_0(\Lambda)$).
Now fix  $\al_0\in\Lambda$ and assume the condition $\al_1=\cdots=\al_m=\bt_m=0$. 
A similar argument as above then gives that $D_n(f(\al_0,\bt_0;z))$ is nonzero
and the identity (\ref{asidbt}) with $\nu=0$ is true for all $\bt_0\in\Lambda$ if $n$ is 
sufficiently large. Continuing this way, we complete the proof of Proposition \ref{asid}
by induction. 

\begin{remark}
A similar argument applies to the asymptotic form of the differential identity (\ref{diffid2new})
we need in Section \ref{secV} below. We omit the discussion.
\end{remark}

\end{proof}

\medskip

We will now complete the proof of Theorem \ref{asTop} in the case $|||\bt|||<1/2$, $V(z)\equiv 0$
by integrating the identities of Proposition \ref{asid}. In this case we denote
\[
D_n(f(z))=D_n(\al_0,\dots,\al_m;\bt_0,\dots,\bt_m).
\]
First, set $m=0$ and $\bt_0=0$. Then (\ref{asidal}) becomes
\be
{\partial\over\partial\al_0}\ln D_n(\al_0)=
2\al_0\left(1+\ln n +{d\over d\al_0}
{\Ga(1+\al_0)\over \Ga(1+2\al_0)}\right)+O(n^{-1}\ln n)+O(\de n\ln n).
\ee
Integrating both sides over $\al_0$ from 0 to some $\al_0$ and using the fact
that $D_n(0)=1$, we obtain
\be\la{Dal0}
D_n(\al_0)=n^{\al_0^2}{G(1+\al_0)^2\over G(1+2\al_0)},
\ee
where $G(x)$ is Barnes $G$-function.
To perform the integration we used the identity
\be\la{idG}
\int_0^z \left(1+{d\over dx}
{\Ga(1+x)\over \Ga(1+2x)}\right)2xdx=\ln{G(1+z)^2\over G(1+2z)},
\ee
which easily follows from the standard formula (see, e.g. \ci{WW}):
\be\la{idst}
\int_0^z\ln \Gamma(x+1)dx=
{z\over 2}\ln 2\pi -{z(z+1)\over2}+z\ln\Gamma(z+1)-\ln G(z+1).
\ee
Now set $m=1$, $\al_0$ fixed. Set $\bt_0=\bt_1=0$.
Relation (\ref{asidal}) for $\nu=1$ is then
\begin{multline}
{\partial\over\partial\al_1}\ln D_n(\al_0,\al_1)=
2\al_1\left(1+\ln n +{d\over d\al_0}
{\Ga(1+\al_1)\over \Ga(1+2\al_1)}\right)+\\
\al_0\ln z_0+\al_0(\ln z_1+i\pi)-2\al_0\ln(z_0-z_1)+2\pi i\al_0+
O(n^{-1}\ln n)+O(\de n\ln n).
\end{multline}
Integrating this over $\al_1$ from 0 to some fixed $\al_1$ along a path lying in $\Lambda$ (see Proposition 
\ref{asid}) and
using (\ref{idG}), we obtain
\begin{multline}
\ln{D_n(\al_0,\al_1)\over D_n(\al_0,0)}=
\al_1^2\ln n +2\ln G(1+\al_1)-\ln G(1+2\al_1)+
\al_0\al_1\ln(z_0z_1e^{3\pi i})\\
-2\al_0\al_1\ln(z_0-z_1)+
O(n^{-1}\ln n)+O(\de n\ln n).
\end{multline}
Substituting here (\ref{Dal0}), we obtain
\begin{multline}
D_n(\al_0,\al_1)=
n^{\al_0^2+\al_1^2}\prod_{j=0}^1{G(1+\al_j)^2\over G(1+2\al_j)}
\left[{(z_0-z_1)^2\over z_0z_1e^{3\pi i}}\right]^{-\al_0\al_1}
(1+O(n^{-1}\ln n)+O(\de n\ln n))
=\\
n^{\al_0^2+\al_1^2}\prod_{j=0}^1{G(1+\al_j)^2\over G(1+2\al_j)}
|z_0-z_1|^{-\al_0\al_1}
(1+O(n^{-1}\ln n)+O(\de n\ln n))
\end{multline}
(to write the last equation we recall (\ref{abszzj}), the way
the branch of $(z-z_j)^{\al_j}$ was fixed there, and the fact that 
$\arg z_0<\arg z_1$). 

Continuing this way, we finally obtain by induction for any fixed $m$,
$\bt_j=0$, $j=0,1,\dots,m$, the asymptotic expression
\be\la{asDal}
D_n(\al_0,\dots,\al_m)=
n^{\sum_{j=0}^m\al_j^2}\prod_{j=0}^m
{G(1+\al_j)^2\over G(1+2\al_j)}
\prod_{0\le j< k\le m}|z_j-z_k|^{-\al_j\al_k}(1+O(n^{-1}\ln n)+O(\de n\ln n)).
\ee

We now add in the $\bt$-singularities. 
We will make use of one more identity, which follows from (\ref{idst}):
\begin{multline}\la{idG2}
\int_0^\bt \left((\al+x){d\over dx}
\ln\Ga(1+\al+x)+(\al-x){d\over dx}\ln\Ga(1+\al-x)-2x
\right)dx\\
=\ln{G(1+\al+\bt)G(1+\al-\bt)\over G(1+\al)^2}.
\end{multline}
First, we obtain the result for the case when $-1/4<\Re\bt_j<1/4$.
This implies that the order of the error term $O(\de n\ln n)$ remains $o(1)$
as we integrate over $\bt_j$'s starting at zero. As before, we always assume integration along 
a path in $\Lambda$. 
Setting $\nu=0$, $\bt_j=0$, $j=1,\dots,m$ in (\ref{asidbt}), and integrating
this identity over $\bt_0$ from zero to a fixed $\bt_0$, we obtain using (\ref{idG2}):
\begin{multline}
\ln{D_n(\al_0,\dots,\al_m;\bt_0)\over D_n(\al_0,\dots,\al_m;0)}=
-\bt_0^2\ln n +\ln{G(1+\al_0+\bt_0)G(1+\al_0-\bt_0)\over G(1+\al_0)^2}\\
+\bt_0\sum_{j\neq 0}\al_j\ln{z_0 e^{i\pi}\over z_j}+
O(n^{-1}\ln n)+O(\de n\ln n).
\end{multline}
Substituting here (\ref{asDal}), we obtain
\begin{multline}\la{542}
D_n(\al_0,\dots,\al_m;\bt_0)=
n^{\sum_{j=0}^m\al_j^2-\bt_0^2}
{G(1+\al_0+\bt_0)G(1+\al_0-\bt_0)\over G(1+\al_0)^2}
\prod_{j=1}^m
{G(1+\al_j)^2\over G(1+2\al_j)}\\
\times\prod_{0\le j< k\le m}|z_j-z_k|^{-\al_j\al_k}
\prod_{j=1}^m \left({z_0 e^{i\pi}\over z_j}\right)^{\al_j\bt_0}
(1+O(n^{-1}\ln n)+O(\de n\ln n)).
\end{multline}
Next, set $\nu=1$, $\bt_j=0$, $j=2,\dots n$ in (\ref{asidbt}) and integrate
over $\bt_1$. We obtain then the determinant
$D_n(\al_0,\dots,\al_m;\bt_0,\bt_1)$. Continuing this procedure,
we finally obtain by induction at step $m$ the asymptotics (\ref{asD})
for the case $-1/4<\Re\bt_j<1/4$
with $V\equiv 0$ and the error term $O(n^{-1}\ln n)+O(\de n\ln n)=o(1)$.

Consider now the general case $|||\bt|||<1/2$. We can choose $q$ such
that  all $\Re\bt_j\in(q-1/4,q+1/4)$.
Divide $(0,q)$ into subintervals of length less than $1/2$.
Apply the above integration procedure to move all $\bt_j$ from zero to the line where the 
real part of all $\bt_j$ is the right end of the first subinterval. Since the length
of subintervals is less than $1/2$, the error term $O(\de\ln n)$ remains $o(1)$. 
Recall that during the integration we avoid any points where $\al_j+\bt_j$ or $\al_j-\bt_j$
is a negative integer.
Next, move the $\bt_j$'s to the right end 
of the second subinterval, and so on, until the point $\Re\bt_j=q$. From that
point move the $\bt_j$'s as needed. We thus obtain  Theorem \ref{asTop} with 
$V\equiv 0$ and $|||\bt|||<1/2$.

\subsection{Pure Fisher-Hartwig singularities. Extension to $|||\bt|||<1$.}\la{pFHext}
We now show that in fact the error term in (\ref{asD}) remains $o(1)$ for
the full range $|||\bt|||<1$. 
First, recall the definition of $\wt R$ and $\om$ in (\ref{Rtilde}) and (\ref{omega}).
We have $|||\bt|||=2\max_j (\Re\bt_j-\om)$ and 
\[
-{1\over 2}<\Re\bt_j-\om<{1\over 2}
\]
for all singular points $z_j$. 
We denote $p_j=z_j$ if $\Re\bt_j-\om>0$, and  $m_+$ the number of such points.
Furthermore,
denote $q_j=z_j$ if $\Re\bt_j-\om<0$, and $m_-$ the number of such points.
Finally, let $r_j=z_j$ if $\Re\bt_j-\om=0$.

Separating the main contributions in $n$ (see (\ref{De1})),
we write the jump matrix for $\wt R$ on $\partial U_{z_j}$ (cf. (\ref{Deas2})) in the form 
\be
I+n^{\om\si_3}\De_1(z)n^{-\om\si_3}+\cdots=
I+\wh\De_1(z)+\wh D(z)+O(n^{-1-\rho}),\qquad
z\in\partial U_{z_j},
\ee
where
\be\la{rho}
\rho=1-|||\bt|||
\ee
and
\be\la{whD}
\wh D(z)=n^{\om\si_3}\De_1(z)n^{-\om\si_3}-\wh\De_1,\qquad z\in\partial U_{z_j},
\ee
\begin{multline}
\wh\De_1(z)=(n^{\om\si_3}\De_1(z)n^{-\om\si_3})_{21}\si_-=
\frac{b_j(z)}{z-p_j}\si_-,\qquad\si_-=\begin{pmatrix} 0 & 0 \cr 1 & 0 \end{pmatrix},\\
b_j(z)= 
-{\Ga(1+\al_j-\bt_j)\over \Ga(\al_j+\bt_j)}
\left({\mathcal{D}(z)\over\ze^{\bt_j}F_j(z)}\right)^{-2}
p_j^{-n} e^{-i\pi(2\bt_j-\al_j)}\frac{z-p_j}{n\ln(z/p_j)}n^{-2\om},
 \qquad z\in\partial U_{p_j},
\end{multline}
\begin{multline}
\wh\De_1(z)=(n^{\om\si_3}\De_1(z)n^{-\om\si_3})_{12}\si_+=
\frac{a_j(z)}{z-q_j}\si_+,\qquad\si_+=\begin{pmatrix} 0 & 1 \cr 0 & 0 \end{pmatrix},\\
a_j(z)= {\Ga(1+\al_j+\bt_j)\over \Ga(\al_j-\bt_j)}
\left({\mathcal{D}(z)\over\ze^{\bt_j}F_j(z)}\right)^2
q_j^n e^{i\pi(2\bt_j-\al_j)}\frac{z-q_j}{n\ln(z/q_j)}n^{2\om},
 \qquad z\in\partial U_{q_j}.
\end{multline}
\be
\wh\De_1(z)=0,\qquad z\in\partial U_{r_j}.
\ee

Note that $\wh D(z)$ and  $\wh\De_1(z)$ are meromorphic functions in 
a neighborhood of $U_{z_j}$ with a simple pole at $z=z_j$.
We have (see (\ref{mub}), (\ref{mu}), and recall that $V\equiv 0$)
\begin{align}\la{b}
b_j=b^{(n)}_j &\equiv\lim_{z\to p_j}b_j(z)=
-n^{2(\bt_j-\om)-1}p_j^{-n+1}\mu_j^{-2}{\Ga(1+\al_j-\bt_j)\over \Ga(\al_j+\bt_j)}=O(n^{-\rho}),\\
a_j=a^{(n)}_j &\equiv\lim_{z\to q_j}a_j(z)=
n^{-2(\bt_j-\om)-1}q_j^{n+1}\mu_j^2{\Ga(1+\al_j+\bt_j)\over \Ga(\al_j-\bt_j)}=O(n^{-\rho}).
\la{a}
\end{align}
Also note that 
\be
\wh D(z)=O(n^{-1}),\qquad z\in\partial U_{z_j}.
\ee

The main idea which will allow us to give the required estimate for the error term in (\ref{asD})
is the following. Write $\wt R$ in the form
\[
\wt R(z)=Q(z)\wh R(z),
\]
where $\wh R(z)$ is the solution to the RHP:
\begin{align}
&\wh R(z)\mbox{ is analytic for } z\in \bbc  \setminus \cup_j\partial U_{z_j}\\
&\wh R(z)_+=\wh R(z)_-(I+\wh\De_1),\qquad z\in\cup_j\partial U_{z_j},\la{whRjump}\\
&\wh R(z)=I+O(1/z),\qquad z\to\infty.
\end{align}

We solve this RHP below explicitly, however, note first that the solution exists and is unique and 
\be\la{Rest}
\wh R(z)=I+O(n^{-\rho})
\ee
uniformly in $z$
by standard arguments. We then have on $\partial U_{z_j}$
using (\ref{Rest}), the jump condition (\ref{whRjump}), and the nilpotency of $\wh\De_1$:
\[
Q_+=\wt R_{+}\wh R_+^{-1}=
\wt R_{-}(I+\wh\De_1(z)+\wh D(z)+O(n^{-1-\rho}))\wh R_+^{-1}=
Q_-(I+\wh D(z)+O(n^{-1-\rho})).
\]
The jump matrix for $Q$ on $\Si^{\mathrm{out}}$, $\Si^{''\mathrm{out}}$ remains exponentially close to the identity. Therefore,
\be\la{Q}
Q(z)=I+Q_1(z)+O(n^{-1-\rho}),\qquad
Q_1(z)={1\over 2\pi i}\int_{\cup_j\partial U_{z_j}}\frac{\wh D(s)}{s-z}ds
\ee

In what follows, we will be interested in the matrix element $Y_{21}(z)$ of our original RHP in a 
neighborhood of $z=0$. In this neighborhood, we have using (\ref{Q}),
\begin{multline}
Y(z)=R(z)N(z)= n^{-\om\si_3}Q(z)\wh R(z)  n^{\om\si_3}N(z)=\\
 n^{-\om\si_3}(I+Q_1(z)+O(n^{-1-\rho}))\wh R(z)  n^{\om\si_3}
\mathcal{D}(z)^{\si_3}
\begin{pmatrix} 0 & 1 \cr -1 & 0 \end{pmatrix}.
\end{multline}
From here, noting that $\mathcal{D}(0)=1$ as $V(z)\equiv 0$, and using the estimates
$\wh R(z)=I+O(n^{-\rho})$, $Q_1(0)=O(n^{-1})$,
we obtain
\[
\chi_{n-1}^{2}=-Y_{21}(0)= 
(1+Q_{1,22}(0)+O(n^{-1-\rho}))\wh R_{22}(0).
\]
Using the expression for $Q_1$ from (\ref{Q}) and  (\ref{whD}), we can write this equation in the form
(recall that the contours $\partial U_{z_j}$ are oriented in the negative direction):
\be\la{chi-new}
\chi_{n-1}^{2}=\left(1-{1\over n}\sum_{j=0}^m(\al_j^2-\bt_j^2)+O(n^{-1-\rho})\right)
\wh R_{22}(0).
\ee
Eventually, we will use the product of these quantities over $n$ to represent the determinant.
Before doing that, we now solve the RHP for $\wh R$ and find $\wh R_{22}(0)$.

Set
\be
\Phi(z)=
\begin{cases}
\wh R(z) & z\in \bbc\setminus \cup_j U_{z_j}\cr
\wh R(z)(I+\wh\De_1(z)) & z\in \cup_j U_{z_j}.
\end{cases}
\ee

So defined, $\Phi(z)$ is obviously a meromorphic function with simple poles at the points $z_j$,
which tends to $I$ at infinity. Therefore, it can be written in the form
\be\la{Phirepr}
\Phi(z)=I+\sum_{j=1}^{m_+}\frac{\Phi_j^+}{z-p_j}+
\sum_{j=1}^{m_-}\frac{\Phi_j^-}{z-q_j}
\ee
for some constant matrices $\Phi^{\pm}_j$. 

Moreover, the function  $\wh R(z)$ (we use here the definition of $\wh\De_1$) in $U_{p_k}$, 
$k=1,\dots, m_+$,
\be
\wh R(z)=\Phi(z)(I+\wh\De_1)^{-1}=
\left[I+\sum_{j=1}^{m_+}\frac{\Phi_j^+}{z-p_j}+
\sum_{j=1}^{m_-}\frac{\Phi_j^-}{z-q_j}\right]\left(I-\frac{b_k(z)}{z-p_k}\si_-\right)
\ee
is analytic in $U_{p_k}$. Hence, the coefficients at negative powers of $z-p_k$ vanish.
Equating the coefficient at $(z-p_k)^{-2}$ to zero, we obtain
\be
\Phi^+_k\si_-=0,
\ee
and therefore the matrix $\Phi^+_k$ has the form
\be\la{b1cond}
\Phi^+_k=\begin{pmatrix} g_k & 0 \cr f_k & 0 \end{pmatrix}
\ee
for some constants $g_k$, $f_k$.
The vanishing of the coefficient at $(z-p_k)^{-1}$ gives the following condition on $\Phi^{\pm}_j$,
where we used the equation $\Phi^+_j\si_-=0$,
\be\la{b2cond}
\Phi_k^+ -b_k\sum_{j=1}^{m_-}\frac{\Phi_j^-\si_-}{p_k-q_j}=b_k\si_-,\qquad k=1,\dots,m_+,
\ee
with $b_k$ given by (\ref{b}). 
Similarly, using the analyticity of $\Phi(z)(I+\wh\De_1(z))^{-1}$ in 
$U_{q_k}$, we obtain for all $k=1,\dots,m_-$
\begin{align}\la{a1cond}
&\Phi^-_k=\begin{pmatrix} 0 & e_k \cr 0 & h_k \end{pmatrix},\\
&\Phi_k^- -a_k\sum_{j=1}^{m_+}\frac{\Phi_j^+\si_+}{q_k-p_j}=a_k\si_+\la{a2cond}
\end{align}
for some constants $e_k$, $h_k$. Conditions (\ref{a2cond}), (\ref{b2cond}) are equations for the 
constants $g_k$, $f_k$, $e_k$, $h_k$. In view of (\ref{chi-new}), we are interested in
\be
\wh R_{22}(0)=\Phi_{22}(0)=
1-\sum_{j=1}^{m_-}\frac{h_j}{q_j},
\ee
where we used (\ref{Phirepr}) to write the last equation.
To calculate this quantity we first substitute $\Phi^+_j$ from (\ref{b2cond})
into (\ref{a2cond}). Then the 
$22$ element of the resulting equations for $k=1,\dots,m_-$ can be written as follows:
\be\la{AB}
(I-A)h=B,
\ee
where $h$ and $B$ are $m_-$-dimensional vectors with components $h_k$ and 
\[
B_k=\sum_{j=1}^{m_+}\frac{a_k b_j}{q_k-p_j},\quad k=1,\dots,m_-,
\]
$A=A^{(n)}$ is an $m_-\times m_-$ matrix with matrix elements
\[
A_{k,\ell}=\sum_{j=1}^{m_+}\frac{a_k b_j}{(q_k-p_j)(p_j-q_{\ell})},
\]
and $I$ is the $m_-\times m_-$ identity matrix. 

Define the  $m_-\times m_-$ diagonal matrix $\De$ as follows 
\[
\De=\mbox{diag}\{-q_1,-q_2,\dots,-q_{m_-}\}.
\]
Then (\ref{AB}) can be written in the form
\[
Tx=y,\qquad T=I-\De^{-1}A\De,\qquad x=\De^{-1}h,\qquad y= \De^{-1}B
\]
By Cramer's rule
\[
x_k=\frac{\det(T_1\cdots y\cdots T_{m_-})}{\det T},
\]
where $T_j$ are the columns of $T$ and $y$ is in the place of the $k$'th column.
We are interested in
\be\la{Tint}
1-\sum_{j=1}^{m_-}\frac{h_j}{q_j}=
1+\sum_{j=1}^{m_-}x_k=\frac{\det(T_1+y,T_2+y,\cdots,T_{m_-}+y)}{\det T}.
\ee
First note that 
\be
\det T= \det (1-\De^{-1}A\De)=\det (I-A)=\det(I-A^{(n)}).
\ee
Second, a direct calculation shows that
\[
(T_1+y,T_2+y,\cdots, T_{m_-}+y)=I-A',\qquad 
A'_{jk}= \sum_{\ell=1}^{m_+}\frac{a'_j b'_{\ell}}{(q_j-p_{\ell})(p_{\ell}-q_k)},
\]
where
\[
a'_j=a_j q^{-1}_j,\qquad b'_{\ell}=b_{\ell}p_{\ell}.
\]
Using the definitions (\ref{a},\ref{b}) of $a_j$, $b_j$, we note that
\[
a'_j=a^{(n-1)}_j+O(n^{-1-\rho}),\qquad b'_j=b^{(n-1)}_j+O(n^{-1-\rho}),
\]
and therefore $A'=A^{(n-1)}+O(n^{-1-\rho})$.
Thus we can rewrite (\ref{Tint})  as
\be\la{Rfinal}
\wh R_{22}(0)=\Phi_{22}(0)=
1-\sum_{j=1}^{m_-}\frac{h_j}{q_j}=\frac{\det(I-A^{(n-1)})}{\det(I-A^{(n)})}
\left[1+O\left({1\over n^{1+\rho}}\right)\right].
\ee

Note that
\[
\det(I-A^{(n)})=1+O(n^{-2\rho}).
\]

Recalling now (\ref{chi-new}) and the representation of the Toeplitz determinant 
$D_n(f)$ as a product of $\chi_k^{-2}$, we can, for some sufficiently large 
$n_0>0$, using (\ref{Rfinal}), write
\begin{multline}\la{newestD}
D_n(f)=D_{n_0}(f)\prod_{k=n_0+1}^{n} \chi_{k-1}^{-2}=\\
 D_{n_0}(f)\prod_{k=n_0+1}^{n}\left[1+{1\over k}\sum_{j=0}^m(\al_j^2-\bt_j^2)
+O\left({1\over k^{1+\rho}}\right)
\right]
\frac{\det(I-A^{(k)})}{\det(I-A^{(k-1)})}
\left[1+O\left({1\over k^{1+\rho}}\right)\right]=\\
C(n_0,\al_0,\dots,\al_m,\bt_0,\dots,\bt_m)n^{\sum_{j=0}^m(\al_j^2-\bt_j^2)}\left[1+
O\left({1\over n^{\rho}}\right)\right],
\end{multline}
where all the error terms are uniform for $\al_j$, $\bt_j$ in compact sets, and 
$C$ is a constant depending analytically on $\al_j$, $\bt_j$, and $n_0$ only.
Recall that in the derivation of this expression we assumed that $|||\bt|||<1$ (and as usual
$\Re\al_j>-1/2$ for all $j$).  Under this condition $\rho>0$ (see (\ref{rho})),
and the error term tends to zero as $n\to\infty$. 
In the previous section, we obtained $C$ explicitly for $\bt_j$ satisfying $|||\bt|||<1/2$.
Obviously, if all $\al_j$, $\bt_j$ belong to fixed compact sets and the value of $n_0$ is fixed, 
the  constant $C$ in (\ref{newestD}) is bounded above in absolute value by a constant 
independent of any $\al_j$, $\bt_j$. By Vitali's theorem this fact implies that
$C$ can be analytically continued in $\bt_j$ to the full domain $|||\bt|||<1$
off its values on the domain $|||\bt|||<1/2$.  Thus $C$ is given by the same expression
also for  $|||\bt|||<1$, and this concludes the proof of Theorem \ref{asTop} for
$|||\bt|||<1$, $V(z)\equiv 0$, with the error term $o(1)=O(n^{|||\bt|||-1})$ in (\ref{asD}).

\subsection{Adding special analytic $V(z)$.}\la{secV}
In this section and in the next one, we will add the multiplicative factor $e^{V(z)}$, where 
$V$ is analytic in a neighborhood of the unit circle $C$, to a
symbol with pure Fisher-Hartwig singularities and obtain the asymptotics of the corresponding determinant.
Consider the deformation of the symbol $f(z,t)$ given by (\ref{ft}) for $t\in[0,1]$.
The analysis is based on integration of the differential identity (\ref{diffid2new}) over $t$.
In the present section we assume that $V$ is such that
\be\la{restr}
1-t+te^{V(z)}\neq 0,\qquad t\in[0,1],\quad z\in C,
\ee
and $1-t+te^{V(z)}$ has no winding around $C$ for all $t\in [0,1]$.
Then the Riemann-Hilbert problem (for the polynomials orthogonal) with $f(z,t)$ has the same singularities
as the problem with $f(z)$ and is solved in the same way. In the following section we remove the condition
(\ref{restr}).

Let us rewrite the identity (\ref{diffid2new}) in terms
of the function $S(z) \equiv S(z,t)$ which is the solution
of the Riemann-Hilbert problem posed on the
deformed contour depicted in Figure 1. From
the transformation $Y \to S$ defined in (\ref{TY}), (\ref{defS})
we obtain
\be \label{Yto S}
Y_{11} = z^{n}f^{-1}S_{12,+} + S_{11,+},\qquad
Y_{21} = z^{n}f^{-1}S_{22,+} + S_{21,+},\qquad z\in C\equiv\Si'.
\ee 
Substituting these expressions into (\ref{diffid2new}) and
taking into account that $\det S=1$ we arrive at the formula:
\begin{multline} \label{differ1000}
{\partial \over \partial t}\ln D_n(f(z,t))=
n\int_{C}f^{-1}\dot{f}\frac{dz}{2\pi iz}
+ \int_{C}\left[-f^{-2}f' +2f^{-1}\left(S'_{22,+}S_{11,+} - S'_{12,+}S_{21,+}\right)
 \right]\dot{f}\frac{dz}{2\pi i}\\
+ \int_{C}\left[z^{-n}\left(S'_{21,+}S_{11,+} - S'_{11,+}S_{21,+} \right)
+ z^{n}\left(S'_{22,+}S_{12,+} - S'_{12,+}S_{22,+} \right)f^{-2}\right] 
\dot{f}\frac{dz}{2\pi i},
\end{multline}
where we introduced the notation
$\dot{f} \equiv \partial f/\partial t$ and   $ f' \equiv \partial f/\partial z$. 
Using the jump relation (\ref{Sjump1}) satisfied by
the function $S(z,t)$ across the unit circle, we can rewrite equation
(\ref{differ1000}) in a more symmetric way, 
\begin{multline} \label{differ100}
{\partial \over \partial t}\ln D_n(f(z,t))=
n\int_{C}f^{-1}\dot{f}\frac{dz}{2\pi iz}
+ X(t),\\
X(t)=
\int_{C}\left[S'_{22,+}S_{11,+} - S'_{12,+}S_{21,+} +
S'_{11,-}S_{22,-} - S'_{21,-}S_{12,-}
\right] f^{-1}\dot{f}\frac{dz}{2\pi i}\\
+ \int_{C}\left[z^{-n}\left(S'_{22,-}S_{12,-} - S'_{12,-}S_{22,-} \right)
+ z^{n}\left(S'_{22,+}S_{12,+} - S'_{12,+}S_{22,+} \right)\right] f^{-2}
\dot{f}\frac{dz}{2\pi i}.
\end{multline}

Analytically continuing the boundary values of $S(z,t)$ from the $``+``$ side of $C$
to the $``-``$ side of $\Si''$, and from the $``-``$ side of $C$
to the $``+``$ side of $\Si$, we can write the term $X(t)$ in
(\ref{differ100}) in the form
\be\la{sumfirst}
X(t)=
\int_{\Si''}(J_- +z^n I_-)f^{-1}\dot{f}\frac{dz}{2\pi i}+
\int_{\Si}(-J_+ +z^{-n} I_+)f^{-1}\dot{f}\frac{dz}{2\pi i},
\ee
where 
\be\la{defIJ}
I=(S'_{22}S_{12} - S'_{12}S_{22})f^{-1},\qquad 
J=S'_{22}S_{11} - S'_{12}S_{21}. 
\ee

Denote by $\Si_{\ep}$ (resp., $\Si''_{\ep}$) 
the part of $\Si$ (resp., $\Si''$) which lies inside $\cup_{j=0}^m U_{z_j}$.
Consider first
\be
\int_{\Si''_{\ep}} (J_- +z^n I_-)f^{-1}\dot{f}\frac{dz}{2\pi i}.
\ee
Using (\ref{Sjump2}) on $\Si''$ and (\ref{defIJ}), we easily obtain that
$J_-=J_+ -z^n I_+$ and $I_-=I_+$, and therefore,
\be\la{intint}
  \int_{\Si''_{\ep}} (J_- +z^n I_-)f^{-1}\dot{f}\frac{dz}{2\pi i}=
  \int_{\Si''_{\ep}} J_ + f^{-1}\dot{f}\frac{dz}{2\pi i}=
  \int_{C''_{\ep}}J f^{-1}\dot{f}\frac{dz}{2\pi i},
\ee
where $C''_{\ep}$ is the part of $\cup_{j=0}^m\partial U_{z_j}$ lying inside the unit circle 
from the intersection of $\partial U_{z_j}$ with the incoming $\Si''$ to the intersection 
with the outgoing $\Si''$ for each $j$. Note that $J(z,t)$ has no nonintegrable singularity at $z_j$.
Indeed, for $z$ inside the smaller sector formed by $\Si''$ at $z_j$, we can write
\be
S(z,t)=Y(z,t)=\wh Y(z,t)
\begin{pmatrix}
1&\kappa(z,t)\cr 0& 1
\end{pmatrix},
\ee
where $\kappa(z,t)=c_1(z,t)(z-z_j)^{2\al_j}$, if $\al_j\neq 0$, and 
$\kappa(z,t)=c_2(z,t)\ln(z-z_j)$, if $\al_j=0$, $\bt_j\neq 0$, for some $c_k(z,t)$ 
analytic near $z_j$, 
chosen so that $\wh Y(z)$ is analytic in a neighborhood of $z_j$. 
Writing $J$ in terms of the matrix 
elements of $\wh Y$, we see that the contributions of the (singular) derivative $\kappa'(z)$
cancel, and we obtain
\be
J=(\wh Y'_{21}\wh Y_{11}-\wh Y_{21}\wh Y'_{11})\kappa+\wh Y'_{22}\wh Y_{11}-\wh Y'_{12}\wh Y_{21}.
\ee

We now analyze (\ref{intint}) asymptotically. 
The asymptotic expression for $S(z,t)$ inside the unit circle and outside $\cup_{j=0}^m U_{z_j}$
is given by (see (\ref{wtR}), (\ref{Ndef}))
\be\la{Sas}
S(z,t)=R(z,t)e^{g(z,t)\si_3}\begin{pmatrix} 0&1\cr -1& 0 \end{pmatrix},
\ee
where $g(z)$ is defined by the formula
\be\la{gg}
\mathcal{D}(z,t)=e^{g(z,t)}.
\ee 
We have
\[
J=g'+R'_{11}R_{22}-R_{12}R'_{21}=g'+O_{\ep}(1/n),\qquad n\to\infty,
\]
and we finally obtain 
\be\la{i1}
\int_{\Si''_{\ep}} (J_- +z^n I_-)f^{-1}\dot{f}\frac{dz}{2\pi i}=
\int_{C''_{\ep}}J f^{-1}\dot{f}\frac{dz}{2\pi i}=
 \int_{C''_{\ep}}g'(z) f^{-1}\dot{f}\frac{dz}{2\pi i}+O_{\ep}(1/n),\qquad n\to\infty.
\ee

Similarly, using the jump condition for $S$ and then the asymptotics for $S$ outside 
the unit circle, we obtain
\be\la{i2}
  \int_{\Si_{\ep}} (-J_+ +z^{-n} I_+)f^{-1}\dot{f}\frac{dz}{2\pi i}=
  -\int_{C_{\ep}}J f^{-1}\dot{f}\frac{dz}{2\pi i}=
  \int_{C_{\ep}}g'(z) f^{-1}\dot{f}\frac{dz}{2\pi i}+O_{\ep}(1/n),\qquad n\to\infty,
\ee
where $C_{\ep}$ is the part of $\cup_{j=0}^m\partial U_{z_j}$ lying outside the unit circle 
from the intersection of $\partial U_{z_j}$ with the incoming $\Si$ to the intersection with 
the outgoing $\Si$ for each $j$. 
 
Returning to the integrals (\ref{sumfirst}), we now consider the part arising from the 
integration over $\Si''\setminus\Si''_{\ep}$ and $\Si\setminus\Si_{\ep}$. In this part 
the terms containing $z^{\pm n}I_{\mp}$ give a contribution which is 
exponentially small in $n$,
while the integration of the terms with $J$ can be replaced by the integration over
$\Si' \setminus\Si'_{\ep}$, where $\Si'_{\ep}=\Si'\cap(\cup_{j=0}^m U_{z_j})$, and over parts 
of the boundaries $\partial U_{z_j}$. Thus, recalling also (\ref{i1}), (\ref{i2}), we have 
\be\la{sumsecond}
X(t)=
\int_{\Si' \setminus\Si'_{\ep}} (J_+ - J_-)f^{-1}\dot{f}\frac{dz}{2\pi i}\\+
\sum_{j=0}^m \left(\int_{\partial U_{z_j}^+}+\int_{\partial U_{z_j}^-}\right)
g'(z) f^{-1}\dot{f}\frac{dz}{2\pi i}+O_{\ep}(1/n),\qquad n\to\infty,
\ee
where $\partial U_{z_j}^+$ (resp., $\partial U_{z_j}^-$)
is the part of the boundary of $U_{z_j}$ inside (resp., outside) the unit circle
oriented from the intersection with the incoming $\Si'$ to the intersection with the outgoing
$\Si'$. 

Note that using the same considerations as before, we can write in (\ref{sumsecond})
$J_+-J_-=g'_+ + g'_-+O_{\ep}(1/n)$, and therefore 
\be\la{sumthird}
X(t)=
\int_{\Si_+}g'_+f^{-1}\dot{f}\frac{dz}{2\pi i}+
\int_{\Si_-}g'_-f^{-1}\dot{f}\frac{dz}{2\pi i}+O_{\ep}(1/n),\qquad n\to\infty,
\ee
where the closed anticlockwise oriented contours 
\be
\Si_+=(\Si'\setminus\Si'_{\ep})\cup_{j=0}^m \partial U_{z_j}^+,\qquad
\Si_-=(\Si'\setminus\Si'_{\ep})\cup_{j=0}^m \partial U_{z_j}^-.
\ee
(Note that one can deform $\Si_+$ (resp., $\Si_-$) to a circle around zero of radius
$1-\ep$ (resp., $1+\ep$).)  

By (\ref{ft}), 
\be\la{ff}
f^{-1}\dot{f} = \frac{-1 + e^{V(z)}}{1-t+te^{V(z)}} = {\partial \over \partial t}
\ln\left(1-t+te^{V(z)}\right).
\ee
Furthermore, writing $g$ in the form
\be\la{ggg}
g(z,t)=g^{Sz}(z,t)+g^{FH}(z),
\ee 
we have, by (\ref{gg}), (\ref{48}), (\ref{Dl1}), (\ref{Dg1}), 
\be\la{gszfh}
g^{Sz}(z,t)=\int_C\frac{\ln(1-t+t e^{V(s)})}{s-z}\frac{ds}{2\pi i},\qquad
g^{FH}(z)=
\begin{cases}
\sum_{k=1}^m (\al_k+\bt_k)\ln\frac{z-z_k}{z_k e^{i\pi}},& |z|<1\\
\sum_{k=1}^m (-\al_k+\bt_k)\ln\frac{z-z_k}{z}, &|z|>1
\end{cases}.
\ee

Note that the solution of the Riemann-Hilbert problem is uniform in $t\in[0,1]$.
From this fact and the explicit formulas (\ref{ff}), (\ref{ggg},\ref{gszfh}), we conclude,
by an argument similar to the argument following (\ref{mark}) that the identity
(\ref{differ100}) holds for all $t\in[0,1]$.
Now integrating equation (\ref{differ100}) from $t =0$ to $t=1$, we connect 
the Toeplitz determinant
$D_n(f(z,1))$ with the Toeplitz determinant $D_{n}(f(z,0))$  which represents
the ``pure'' Fisher-Hartwig case and whose asymptotics  we
evaluated in the previous section. First, using (\ref{ff}) and changing the order 
of integration in the first term of (\ref{differ100}), we obtain
\be\la{int1}
\int_0^1 dt\int_C f^{-1}\dot{f}\frac{dz}{2\pi i z}=
\frac{1}{2\pi}\int_{0}^{2\pi}V(e^{i\theta})d\theta=V_0.
\ee
Furthermore, by (\ref{sumthird}), (\ref{ggg},\ref{gszfh}),
\be\la{int2}
\int_0^1dt X(t)= I^{Sz}+I^{FH}+O_{\ep}(1/n),\qquad n\to\infty,
\ee
where (cf. \ci{DeiftInOp}, Eq (86),(87))
\be\la{isz}
I^{Sz}=\int_0^1 dt\int_C\left((g^{Sz})'_+ + (g^{Sz})'_-\right)f^{-1}\dot{f}\frac{dz}{2\pi i}=
\sum_{k=1}^\infty k V_k V_{-k},
\ee
and
\begin{multline}
I^{FH}=\int_{\Si_+}\left(g^{FH}(z)\right)'_+ V(z)\frac{dz}{2\pi i}+
\int_{\Si_-}\left(g^{FH}(z)\right)'_- V(z) \frac{dz}{2\pi i}\\=
\sum_{k=0}^m \left[
(\al_k+\bt_k)\int_{\Si_+}\frac{V(z)}{z-z_k}\frac{dz}{2\pi i}+
(-\al_k+\bt_k)\int_{\Si_-}\left(\frac{V(z)}{z-z_k}-\frac{V(z)}{z}\right)
\frac{dz}{2\pi i}\right].
\end{multline}

Since
\be
g^{Sz}(z,1)=\int_C\frac{V(s)}{s-z}\frac{ds}{2\pi i}=
\begin{cases}
\ln b_+(z)+V_0,& |z|<1\\
-\ln b_-(z), &|z|>1
\end{cases},
\ee
we obtain
\be
\ln b_+(z_k)=\int_{\Si_-}\frac{V(z)}{z-z_k}\frac{dz}{2\pi i}-V_0,\qquad
\ln b_-(z_k)=-\int_{\Si_+}\frac{V(z)}{z-z_k}\frac{dz}{2\pi i},
\ee
which finally gives
\be\la{ifh}
I^{FH}=
\sum_{k=0}^m \left[
-(\al_k+\bt_k)\ln b_-(z_k)+(-\al_k+\bt_k)\ln b_+(z_k)\right].
\ee

Collecting (\ref{int1}), (\ref{int2}), (\ref{isz}), and (\ref{ifh}), we obtain from (\ref{differ100})
\begin{multline}
\ln D_n(f(z,1)) - \ln D_n(f(z,0)) = nV_{0}+
\sum_{k=1}^\infty k V_k V_{-k}\\+
\sum_{k=0}^m \left[
-(\al_k+\bt_k)\ln b_-(z_k)+(-\al_k+\bt_k)\ln b_+(z_k)\right]
+O_{\ep}(1/n),\qquad n\to\infty,
\end{multline}
which, in view of the result of the previous section, concludes the proof of Theorem \ref{asTop}
for analytic $V(z)$ satisfying the condition (\ref{restr}).

\subsection{Extension to general analytic $V(z)$.}\la{secVGenAn}
Now let $V(z)$ be any function analytic in a neighborhood of the unite circle. 
Since zeros of the expression $1-t+t e^{V(z)}$, $t\in[0,1]$, $z\in C$, can only occur
if $\Im V(z)=\pi (2k+1)$, $k\in\bbz$, 
there exists
a positive integer $q$ such that $\frac{1}{q}V(z)$ satisfies the condition (\ref{restr})
of the previous section, i.e., 
\[
1-t + t e^{\frac{1}{q}V(z)} \neq 0,\quad \forall t \in [0,1], \quad z \in C,
\]
and this function has no winding around $C$ for all $t \in [0,1]$.
Let
\begin{equation}\label{fl}
f_{0}(z) = f^{FH}(z), \quad
f_{\ell}(z) = e^{\frac{1}{q}V(z)}f_{\ell-1}(z), \quad \ell =1,\dots, q,
\end{equation}
where $f^{FH}(z)$ is the symbol for the ``pure'' Fisher-Hartwig case.
Note that 
\begin{equation}\label{fp}
f(z) \equiv  e^{V(z)}f^{FH}(z) = f_{q}(z).
\end{equation}
Consider $f_{\ell}(z)$, $\ell=1,\dots,q$, and introduce the deformation,
\begin{equation}\label{fpt}
f_\ell(z,t)  =  \left(1-t +te^{\frac{1}{q}V(z)}\right)f_{\ell-1}(z)=
\left(1-t +te^{\frac{1}{q}V(z)}\right)e^{\frac{\ell-1}{q}V(z)}f^{FH}(z).
\end{equation} 
All the considerations  of the part of the  previous section
between equations (\ref{Yto S}) and (\ref{sumthird}) go through with $f(z,t)$ replaced by $f_\ell(z,t)$ 
and we arrive at the formulae:
\begin{multline} \label{differ100p}
{\partial \over \partial t}\ln D_n(f_\ell(z,t))=
n\int_{C}f_\ell^{-1}\dot{f_\ell}\frac{dz}{2\pi iz}
+ X(t),\\
X(t)=  
\int_{\Si_+}g'_+f_\ell^{-1}\dot{f_\ell}\frac{dz}{2\pi i}+
\int_{\Si_-}g'_-f_\ell^{-1}\dot{f_\ell}\frac{dz}{2\pi i}+O_{\ep}(1/n),\qquad n\to\infty,
\end{multline}
where, as before, the closed anticlockwise oriented contours 
\be
\Si_+=(\Si'\setminus\Si'_{\ep})\cup_{j=0}^m \partial U_{z_j}^+,\qquad
\Si_-=(\Si'\setminus\Si'_{\ep})\cup_{j=0}^m \partial U_{z_j}^-,
\ee
and $g(z)\equiv g_{\ell}(z)$ now corresponds to $f_{\ell}$.

The first term in (\ref{differ100p}) yields (cf. (\ref{int1}))
\begin{equation}\label{int1p}
\int_0^1 dt\int_C f_\ell^{-1}\dot{f_\ell}\frac{dz}{2\pi i z}=
\frac{1}{2\pi q}\int_{0}^{2\pi}V(e^{i\theta})d\theta=\frac{1}{q}V_0.
\end{equation}
In order to evaluate $X(t)$, we write $g$ in the form
\be\label{gp}
g(z,t)=g^{Sz}(z,t)+ \wt{g}^{Sz}(z) + g^{FH}(z),
\ee 
where $g^{FH}(z)$ is the same as in (\ref{gszfh}), and
\be\label{gzp}
g^{Sz}(z,t)=\int_C\frac{\ln(1-t+t e^{\frac{1}{q}V(s)})}{s-z}\frac{ds}{2\pi i},\qquad
\wt{g}^{Sz}(z)=\frac{\ell-1}{q}\int_C\frac{V(s)}{s-z}\frac{ds}{2\pi i}.
\ee
Then we obtain
\be\la{int2p}
\int_0^1dt X(t)= I^{Sz}+ \wt{I}^{Sz} + I^{FH}+O_{\ep}(1/n),\qquad n\to\infty,
\ee
where, up to the replacement $ V \to \frac{1}{q}V$, the 
integrals  $I^{Sz}$ and  $I^{FH}$ are the respective integrals from the
previous section, i.e.,
\begin{equation}\label{ISzFHp}
I^{Sz}= \frac{1}{q^2}\sum_{k=1}^\infty k V_k V_{-k},\quad
I^{FH}=
\frac{1}{q}\sum_{k=0}^m \left[
-(\al_k+\bt_k)\ln b_-(z_k)+(-\al_k+\bt_k)\ln b_+(z_k)\right].
\end{equation}
The term $\wt{I}^{Sz}$ in (\ref{int2p}) is given by the equation
\be\label{Iszp}
\wt{I}^{Sz} = \frac{1}{q}\int_{C}\Bigl(\left(\wt{g}^{Sz}(z)\right)'_+ + 
\left(\wt{g}^{Sz}(z)\right)'_-\Bigr)V(z)\frac{dz}{2\pi i}.
\ee
Note that
\[
\left(\wt{g}^{Sz}(z)\right)'_+ = \frac{\ell-1}{q}\sum_{k=1}^{\infty}kz^{k-1}V_k,\qquad
\left(\wt{g}^{Sz}(z)\right)'_- = \frac{\ell-1}{q}\sum_{k=1}^{\infty}kz^{-k-1}V_{-k}.
\]
Therefore, after a simple calculation we obtain 
\be\la{Iszp2}
\wt{I}^{Sz} = \frac{2\ell-2}{q^2}\sum_{k=1}^{\infty}kV_{k}V_{-k}.
\ee

Integrating (\ref{differ100p}) from $t =0$ to $t=1$ and taking into account
(\ref{int1p}), (\ref{int2p}),  (\ref{ISzFHp}), and (\ref{Iszp2}), we
obtain the following equation for the determinant $D_n(f_\ell(z))$:
\begin{multline}\label{final0}
\ln D_n(f_\ell(z)) - \ln D_n(f_{\ell-1}(z)) = \frac{1}{q}nV_{0}+
\frac{2\ell-1}{q^2}\sum_{k=1}^\infty k V_k V_{-k}\\+
\frac{1}{q}\sum_{k=0}^m \left[
-(\al_k+\bt_k)\ln b_-(z_k)+(-\al_k+\bt_k)\ln b_+(z_k)\right]
+O_{\ep}(1/n),\qquad n\to\infty.
\end{multline}
This equation holds for any $\ell=1,\dots, q$.
Summing up from $\ell =1$ to $\ell=q$ we again arrive at the formula
\begin{multline}\label{final1}
\ln D_n(f(z)) - \ln D_n(f^{FH}(z)) = nV_{0}+
\sum_{k=1}^\infty k V_k V_{-k}\\+
\sum_{k=0}^m \left[
-(\al_k+\bt_k)\ln b_-(z_k)+(-\al_k+\bt_k)\ln b_+(z_k)\right]
+O_{\ep}(1/n),\qquad n\to\infty,
\end{multline}
which conludes the proof of Theorem \ref{asTop} in the case of $V(z)$ analytic
in a neighborhood of the unit circle.

\subsection{Extension to smooth $V(z)$}\la{secVext}
If $V(z)$ is just sufficiently smooth, in particular $C^\infty$,
 on the unit circle $C$ so that (\ref{Vcond})
holds for $s$ from zero up to and including
some $s\ge 0$, 
we can approximate $V(z)$ by trigonometric polynomials
$V^{(n)}(z)=\sum_{k=-p(n)}^{p(n)} V_k z^k$, $z\in C$. 
First, consider the case when 
$|||\bt|||=\max_{j,k}|\Re\bt_j-\Re\bt_k|=2\max_j |\Re\bt_j-\om|<1$, 
where $\om$ is defined by (\ref{omega}). (The indices $j,k=0$ are omitted if $\al_0=\bt_0=0$.)
We set 
\be\la{pnu}
p=[n^{1-\nu}],\qquad \nu=|||\bt|||+\ep_1,
\ee
where $\ep_1>0$ is chosen sufficiently small so that $\nu<1$
(square brackets denote the integer part).

First, we need to extend the RH analysis of the previous sections to symbols which
depend on $n$, namely to the case when $V$ in $f$ is replaced by $V^{(n)}$. (We will denote
such $f$ by $f(z,V^{(n)})$, and the original $f$ by $f(z,V)$.)
We need to have a suitable estimate for the behavior of the error term in the asymptotics with $n$. 
For a fixed $f$, our analysis depended, in particular, on the fact that $f(z)^{-1}z^{-n}$ 
is of order
$e^{-\ep' n}$, $\ep'>0$, for $z\in \Si^\mathrm{out}$ (see Section \ref{RRHP}), and similarly,
$f(z)^{-1}z^n=O(e^{-\ep' n})$ for  $z\in \Si^{''\mathrm{out}}$.
Here the contours $\Si^\mathrm{out}$, $\Si^{''\mathrm{out}}$ are outside a {\it fixed} 
neighborhood of the unit circle (outside and inside $C$, respectively). 
If $V$ is replaced by $V^{(n)}$, let us define the curve $\Si$ outside $\cup_{j=0}^m U_{z_j}$ by
\be\la{si1}
z=\left(1+\ga\frac{\ln p}{p}\right)e^{i\th},\qquad \ga>0, 
\ee
and $\Si^{''}$ outside $\cup_{j=0}^m U_{z_j}$ by 
\be\la{si2}
z=\left(1-\ga\frac{\ln p}{p}\right)e^{i\th}.
\ee
Inside all the sets $U_{z_j}$, the curves still go to $z_j$ as discussed in Section \ref{RHa}.
Let the radius of all $U_{z_j}$ be $2\ga\ln p/p$.
We now fix the value of $\ga$ as follows.
Using the condition (\ref{Vcond}) we can write (here and below $c$ stands for various 
positive constants independent of $n$)
\begin{multline}
|V^{(n)}(z)|-|V_0|
\le\sum_{k=-p,\;k\neq 0}^p |k^s V_k| \frac{|z|^k}{|k|^s}
< c\left(\sum_{k=-p,\;k\neq 0}^p |k^s V_k|^2\right)^{1/2}
\left(\sum_{k=1}^p\frac{(1\pm 3\ga\ln p/p)^{\pm2k}}{k^{2s}}\right)^{1/2}\\
<c\left(\sum_{k=1}^p\frac{(1\pm 3\ga\ln k/k)^{\pm2k}}{k^{2s}}\right)^{1/2}<
c\left(\sum_{k=1}^p{1\over k^{2(s-3\ga)}}
\left[1+O\left({\ln^2 k\over k}\right)\right]\right)^{1/2},
\end{multline}
where $z\in\Si^\mathrm{out}$, $z\in\partial U_{z_j}\cap\{|z|>1\}$  
(with ``$+$'' sign in ``$\pm$''), and $z\in \Si^{''\mathrm{out}}$,
$z\in\partial U_{z_j}\cap\{|z|<1\}$
(with ``$-$'' sign). We now set
\be\la{gas}
3\ga=s-(1+\ep_2)/2,\qquad \ep_2>0,
\ee
and then
\be\la{Vbound}
|V^{(n)}(z)|<c,\qquad |b_+(z,V^{(n)})|<c,\qquad |b_-(z,V^{(n)})|<c,\qquad
\mbox{for all } n
\ee
uniformly on $\Si^\mathrm{out}$, $\Si^{''\mathrm{out}}$, $\partial U_{z_j}$'s, and
in fact in the whole annulus $1-3\ga\frac{\ln p}{p}<|z|<1+3\ga\frac{\ln p}{p}$.

It is easy to adapt the considerations of the previous sections to the present case, and 
we again obtain the expansion (\ref{Deas}) for the jump matrix of $R$ on $\partial U_{z_j}$.
Note that now $|\ze(z)|=O(n^\nu\ln n)$ and $|z-z_j|=\ln n/ n^{1-\nu}$ as $n\to\infty$
for $z\in\partial U_{z_j}$, and therefore using (\ref{Deas}), (\ref{Fj}), 
(\ref{abszzj}),
(\ref{Dl1}) and the definition of $\nu$ in (\ref{pnu}), we obtain, in particular,
\be
n^{\om\si_3}\De_1(z)n^{-\om\si_3}=O\left({1\over n^{\ep_1}\ln n}\right),\qquad 
z\in\cup_{j=0}^m\partial U_{z_j}.
\ee
Furthermore, as follows from (\ref{si1}), (\ref{si2}), (\ref{Vbound}), 
and (\ref{Rs1},\ref{Rs2}),
the jump matrix on $\Si^\mathrm{out}$ and $\Si^{''\mathrm{out}}$ is now the identity 
plus a function uniformly bounded in absolute value by 
\be\la{bbb}
c\left({n^{1-\nu}\over\ln n}\right)^{2\max_j|\Re\bt_j|}
\left(1\pm\ga(1-\nu)\frac{\ln n}{n^{1-\nu}}\right)^{\mp n}< 
c\exp\left\{-{\ga\over2}(1-\nu)n^\nu\ln n\right\}n^{2(1-\nu)\max_j|\Re\bt_j|},
\ee
where the upper sign corresponds to  $\Si^\mathrm{out}$, and the lower to 
$\Si^{''\mathrm{out}}$.

The RH problem for $R(z)$ (see Section \ref{RRHP}) 
is therefore solvable,  and we obtain $R(z)$ as a series where the first term $R_1$ is the same
as before, and for the error term the same estimate holds
for $z$ outside a fixed neighborhood of the unit circle, e.g., for $z$ large.

This implies that Theorem \ref{asTop} holds for $f(z,V^{(n)})$. Note that it also holds for
$|f(z,V^{(n)})|$.

We will now show that replacing $V^{(n)}$ with $V$ in the symbol of the determinant 
$D_n(f(z,V^{(n)}))$ results, under a condition on $s$, 
in a small error only, so that Theorem \ref{asTop} holds for $D_n(f(z,V))$ as well.

Using the Heine representation (\ref{ir}) for a Toeplitz determinant with (any) symbol $f(z)$,
the straightforward estimate 
\be\la{bb}
b_\pm(z, V^{(n)})=b_\pm(z, V)\left[1+O\left({1\over n^{(1-\nu)s}}\right)\right],
\qquad \mbox{uniformly for}\quad |z|=1,
\ee
which follows from (\ref{Vcond}),
and Theorem \ref{asTop} for $D_n(| f(z,V^{(n)})|)$
and $D_n(f(z,V^{(n)}))$, we have if $s(1-\nu)>1$,
\begin{multline}\la{estD}
\left|D_n(f(z,V))-D_n(f(z,V^{(n)})\right|<\\
\frac{1}{(2\pi)^{n}n!}\int_0^{2\pi}\cdots\int_0^{2\pi}
\prod_{1\leq j < k \leq n}|e^{i\phi_{j}} - e^{i\phi_{k}}|^{2}
\prod_{j=0}^n | f(e^{i\phi_j},V^{(n)})|d\phi_{j}\times
\left(\left|1+c/n^{(1-\nu)s}\right|^n-1\right)\\
<c e^{\Re V_0 n} n^{\sum_{j=0}^m((\Re\al_j)^2+(\Im\bt_j)^2)}
(e^{c/n^{(1-\nu)s-1}}-1)\\
<c \left|e^{V_0 n} n^{\sum_{j=0}^m (\al_j^2-\bt_j^2)}\right|
n^{\sum_{j=0}^m((\Im\al_j)^2+(\Re\bt_j)^2)}
{1\over n^{(1-\nu)s-1}}\\
<c\left|D_n(f(z,V^{(n)})\right|n^{-((1-\nu)s-1-\sum_{j=0}^m((\Im\al_j)^2+(\Re\bt_j)^2))}.
\end{multline}
Therefore,
\be\la{estD2}
D_n(f(z,V))=D_n(f(z,V^{(n)}))\left(1+\frac{D_n(f(z,V))-D_n(f(z,V^{(n)}))}
{D_n(f(z,V^{(n)}))}\right)=D_n(f(z,V^{(n)}))(1+o(1)),
\ee
if
\be\la{s-est}
s>\frac{1+\sum_{j=0}^m((\Im\al_j)^2+(\Re\bt_j)^2)}{1-\nu}.
\ee
Note that this condition is consistent with (\ref{gas}) and the requirement that $\ga>0$.
Using the expression for $\nu$ in (\ref{pnu}) and noting that  $\ep_1$ can be arbitrary close to zero, 
we replace (\ref{s-est}) with  (\ref{s-main0}).
Under the condition (\ref{s-main0})  we then obtain the statement of the theorem for
$D_n(f(z,V))$.

\section{Appendix. The Toeplitz determinant $D_n$ as a tau-function}
In this section we construct a Fuchsian system of ODE's corresponding to the Riemann-Hilbert
problem of Section \ref{RHsection} for $V\equiv 0$. We show that the differential identities 
(\ref{diffid1}) for the Toeplitz determinant can be viewed as monodromy deformations 
of the tau-function associated with this Fuchsian system.

Assume the pure Fisher-Hartwig case, $V(z) \equiv 0$.  Set
\begin{equation}\label{Psitau}
\Phi(z) = \Lambda
Y^{(n)}(z) \Lambda^{-1} \prod_{k=0}^{m}(z-z_k)^{\alpha_k\sigma_3}z^{\lb\sigma_3},
\end{equation}
where
$$
\Lambda =  \prod_{k=0}^{m}z_k^{\frac{\beta_k + \alpha_k}{2}\sigma_3},
\qquad  \lb = \sum_{k=0}^{m}\frac{\beta_k - \alpha_k}{2} - \frac{n}{2},
$$
and the branches of all multi-valued functions are chosen as in Section \ref{diffid}.
In terms of the function $\Phi(z)$, the Riemann-Hilbert problem (\ref{RHPYb})-
(\ref{RHPYc}) reads as follows:
\begin{enumerate}
    \item[(a)]
        $\Phi(z)$ is  analytic for $z\in\bbc \setminus(C \cup [0, 1]\cup \{\cup_{j=0}^{m}\Gamma_j\})$, where
         $\Gamma_j$ is
        the ray $\theta = \theta_j$ from $z_j$ to infinity. The unit circle $C$ is oriented 
        as  before, counterclockwise,  the segment $[0, 1]$ is oriented from
        $0$ to $1$, and the rays $\Gamma_j$ are oriented towards infinity.
    \item[(b)] The boundary values of $\Phi(z)$ are related by the jump conditions,
\begin{equation}\label{RHPPsib}
            \Phi_+(z) = \Phi_-(z)
   \begin{cases}
 \begin{pmatrix}
                1 & s_j \cr
                0 & 1
             \end{pmatrix}, & z\in C, \quad \theta_j < \arg z < \theta_{j+1},\quad 
             0\leq j \leq m, \,\, \theta_{m+1} = 2\pi\cr\cr
            e^{-2\pi i\alpha_j\sigma_3}, & z\in \Gamma_j, \quad 1\leq j \leq m\cr\cr
e^{-2\pi i(\alpha_0 + \lb)\sigma_3}, & z\in \Gamma_0\cr\cr
 e^{-2\pi i\lb\sigma_3}, & z\in [0,1]\cr\cr             
\end{cases}
 \end{equation}
 where
 \be\la{ajdef0}
s_j = \exp\left\{-i\pi\sum_{k=0}^{j}\beta_k + i\pi\sum_{k=j+1}^{m}\beta_k
-i\pi\sum_{k=0}^{j}\alpha_k - 3i\pi\sum_{k=j+1}^{m}\alpha_k\right\}.
\ee
(for $j =m$, the second and the fourth sums are absent)
\item[(c)] 
        $\Phi(z)$ has the following asymptotic behavior at infinity:
        \begin{equation} \label{RHPPsic}
            \Phi(z) = \left(I+ O \left( \frac{1}{z} \right)\right)
            z^{\left(\sum_{k=0}^{m}\frac{\beta_k +\alpha_k}{2}+\frac{n}{2}\right)\sigma_3}
            e^{2\pi i\sum_{k=j+1}^{m}\alpha_k\sigma_3},
                   \end{equation}
as $z \to \infty$  and $\theta_{j} < \arg z < \theta_{j+1}$, $0\leq j\leq m$, $\theta_{m+1} =2\pi$
(for $j =m$ the last factor is omitted).
\item [(d)] In the neighborhoods $U_{z_j}$ of the points $z_j$,
    $j =0, 1,\dots, m$, the function $\Phi(z)$ admits the following
    representations, which constitute a refinement of the estimates (\ref{RHPYd})
    and (\ref{RHPYe}).
 \begin{itemize}
 \item If $\alpha_j \neq 0$, then
    \be\label{Psij0}
\Phi(z) = \wt{\Phi}_{j}(z)(z-z_j)^{\alpha_j\sigma_3}C_j,
\ee
where $\wt{\Phi}_j(z)$ is holomorphic at $z = z_j$ (it is essentially
the function $\wt{Y}(z)$ from Section \ref{diffid})  and the matrix $C_j$ is
given by the formula,
\be\la{Cdef10}
C_j = \begin{pmatrix} 1 & c_j\cr
0&1
\end{pmatrix},
\ee
with
\be\label{Cdef20}
c_j = s_j\begin{cases}
\frac{1 - e^{2\pi i(\beta_j + \alpha_j)}}{1 -e^{4\pi i\alpha_j}}, & z\in U_{z_j},\quad |z| <1\cr\cr
\frac{1 - e^{2\pi i(\beta_j - \alpha_j)}}{1 -e^{4\pi i\alpha_j}}, & z\in U_{z_j},\quad |z| >1, \quad \arg z < \theta_j \cr\cr
e^{4\pi i\alpha_j}\frac{1 - e^{2\pi i(\beta_j - \alpha_j)}}{1 -e^{4\pi i\alpha_j}}, & 
z\in U_{z_j},\quad |z| >1, \quad \arg z > \theta_j
\end{cases}
\ee
in the case $j\neq 0$, and
\be\la{Cdef100}
C_0 = \begin{pmatrix} 1 & c_0\cr
0&1
\end{pmatrix}\times
\begin{cases}
I,&\Im z >0\cr
e^{2\pi i\kappa\sigma_3},&\Im z <0
\end{cases}
\ee
with
\be\label{Cdef50}
c_0 = s_0\begin{cases}
\frac{1 - e^{2\pi i(\beta_0 + \alpha_0)}}{1 -e^{4\pi i\alpha_0}}, & z\in U_{z_0},\quad |z| < 1\cr\cr
\frac{1 - e^{2\pi i(\beta_0 - \alpha_0)}}{1 -e^{4\pi i\alpha_0}}
 &  z\in U_{z_0},\quad |z| >1,\quad \Im z <0 \cr\cr
e^{4\pi i\alpha_0}\frac{1 - e^{2\pi i(\beta_0 - \alpha_0)}}{1 -e^{4\pi i\alpha_0}}, & 
 z\in U_{z_0},\quad |z| >1,\quad \Im z >0
\end{cases}
\ee
in the case $j=0$.
\item If $\alpha_j =0$ and $\beta_j \neq 0$, then
    \be\label{Psij01}
\Phi(z) = \wt{\Phi}_{j}(z)\begin{pmatrix}1&\frac{d_j}{2\pi i}\ln(z-z_j)\cr
0&1\end{pmatrix}C_j,\qquad d_j = -s_j(1 -e^{2\pi i\beta_j}),
\ee
where $\wt{\Phi}_j(z)$ is again holomorphic at $z = z_j$  and the matrix $C_j$
this time is given by the formula
\be\la{Cdef101}
C_j = \begin{pmatrix} 1 & c_j\cr
0&1
\end{pmatrix},
\ee
with
\be\label{Cdef201}
c_j = \begin{cases}
s_{j-1}, & z\in U_{z_j},\quad |z| <1\cr\cr
0 & z\in U_{z_j},\quad |z| >1, \quad \arg z < \theta_j \cr\cr
d_j, & 
z\in U_{z_j},\quad |z| >1, \quad \arg z > \theta_j
\end{cases}
\ee
in the case $j\neq 0$, and
\be\la{Cdef1000}
C_0 = \begin{pmatrix} 1 & c_0\cr
0&1
\end{pmatrix}\times
\begin{cases}
I,&\Im z >0\cr
e^{2\pi i\lb\sigma_3},&\Im z <0
\end{cases}
\ee
with
\be\label{Cdef501}
c_0 = \begin{cases}
s_m e^{4\pi i\lb}, & z\in U_{z_0},\quad |z| < 1\cr\cr
0 &  z\in U_{z_0},\quad |z| >1,\quad \Im z <0 \cr\cr
d_0, & 
 z\in U_{z_0},\quad |z| >1,\quad \Im z >0
\end{cases},
\ee
in the case $j=0$.
\end{itemize} 
\item [(e)] In a small neighborhood of $z =0$, the function $\Phi(z)$
admits a similar representation:
\be\label{zerorep}
\Phi(z) = \wt{\Phi}^{(0)}(z) z^{\lb\sigma_3},
\ee
where $\wt{\Phi}^{(0)}(z)$ is holomorphic at $z=0$.
 \end{enumerate}
We also note that  all the matrices above have determinants equal to $1$.

A key feature of the $\Phi$-RH problem is that all its jump matrices
and the connection matrices $C_j$ are piecewise constant in $z$. 
By standard arguments (see, e.g., \cite{JMU} , \cite{ITW}), based on the Liouville theorem, this fact implies that
the function $\Phi(z)$ satisfies a linear matrix ODE of Fuchsian type:
\be\label{fuchsian}
\frac{d\Phi(z)}{dz} = A(z)\Phi(z), \quad A(z) = \sum_{k=0}^{m}\frac{A_k}{z-z_k} + \frac{B}{z}
\ee
with
\be\label{fuchsian1}
B = \lb\wt{\Phi}^{(0)}(0)\sigma_3\Bigl(\wt{\Phi}^{(0)}(0)\Bigr)^{-1}
\ee
and
\be\label{fuchsian2}
A_j= \begin{cases}\alpha_j\wt{\Phi}_{j}(z_j)\sigma_3\wt{\Phi}^{-1}_{j}(z_j),&\mbox{if}\quad
\alpha_j \neq 0\cr\cr
\frac{d_j}{2\pi i}\wt{\Phi}_{j}(z_j)\begin{pmatrix}0 & 1\cr
0&0\end{pmatrix}\wt{\Phi}^{-1}_{j}(z_j),&\mbox{if}\quad
\alpha_j = 0,\quad \beta_j \neq 0.\end{cases}
\ee
Moreover, since  the jump matrices
and the connection matrices $C_j$ are all constant with respect
to $z_j$, the function $\Phi(z)$  satisfies, in addition to (\ref{fuchsian}),  the equations
\be\label{fuchs3}
\frac{\partial{\Phi(z)}}{\partial z_j} = -\frac{A_j}{z-z_j}\Phi(z),\quad j =1,\dots, m.
\ee
The compatibility condition of (\ref{fuchsian}) and (\ref{fuchs3}) yields the
following nonlinear systems of ODEs on the matrix
coefficients $B$ and $A_j$:
\be\label{Schles1}
\frac{\partial{B}}{\partial z_j}= \frac{[A_j, B]}{z_j},\quad
\frac{\partial{A_k}}{\partial z_j}= \frac{[A_j, A_k]}{z_j-z_k}, \quad k \neq j
\ee
\be\label{Schles2}
\frac{\partial{A_j}}{\partial z_j}= -\sum_{k=0\atop k \neq j}^{m}\frac{[A_j, A_k]}{z_j-z_k}
+ \frac{[B, A_j]}{z_j}.
\ee
In the context of the Fuchsian system (\ref{fuchsian}), the jump matrices and the connection
matrices $C_j$ of the $\Phi$-Riemann-Hilbert
problem form the monodromy data of the system. The fact that these
data do not depend on the parameters $z_j$ means that the functions $B\equiv B(z_1,\dots, z_m$ and 
$A_j \equiv A_{j}(z_1,\dots, z_m)$,
$j =0,\dots, m$ describe {\it isomonodromy  deformations} of  the system  (\ref{fuchsian}). Equations (\ref{Schles1})--(\ref{Schles2}) are the classical Schlesinger equations.

An important role in the modern theory of  isomonodromy
deformations is played by the notion of a $\tau$-function
which  was introduced by M. Jimbo, 
T. Miwa and K. Ueno in \cite{JMU}. In the Fucshian case, the $\tau$-function is defined as follows.
Let
\be\label{tau1}
\wt\omega =  \sum_{k=1}^{m}\mbox{Res}_{z=z_k}\mbox{trace}\,
A(z)\frac{\partial\wt{\Phi}_k(z)}{\partial z}
\wt{\Phi}^{-1}_k(z)dz_k = 
\sum_{k=1}^{m}\mbox{trace}\,
A_k\frac{\partial\wt{\Phi}_k(z_k)}{\partial z}
\wt{\Phi}^{-1}_k(z_k)dz_k.
\ee
As shown in  \cite{JMU}, the differential form 
$\wt\omega \equiv \wt\omega\left(A_0,\dots, A_m, B; z_1,\dots, z_m\right)$
is closed on the solutions of the Schlesinger system (\ref{Schles1})--(\ref{Schles2}).
The $\tau$-function is then defined as the exponential of the antiderivative of $\wt\omega$, i.e.,
\be\label{tau2}
\frac{\partial \ln\tau}{\partial z_k} = \mbox{trace}\,
A_k\frac{\partial\wt{\Phi}_k(z_k)}{\partial z}
\wt{\Phi}^{-1}_k(z_k).
\ee
It was already observed (see, e.g., \cite{ITW} and \cite{BER})  
that the $\tau$-functions evaluated on solutions of the
Schlesinger equations generated by the Riemann-Hilbert problems associated
with Toeplitz and Hankel determinants coincide,  up to trivial factors, with the determinants 
themselves. In particular,  for the Toeplitz determinant
$D_n$ with the pure Fisher-Hartwig symbol,  $V(z) \equiv 0$,  one can follow the
calculations of \cite{ITW} and obtain that
\be\label{tau3}
\frac{\partial \ln D_n}{\partial z_k} =  \mbox{trace}\,
A_k\frac{\partial\wt{\Phi}_k(z_k)}{\partial z}
\wt{\Phi}^{-1}_k(z_k) + 2\sum_{j=0\atop j\neq k}^m\frac{\alpha_k\alpha_j}{z_j-z_k}
-2\lb\frac{\alpha_k}{z_k}.
\ee
Therefore, in the case of the 
Riemann-Hilbert problem (\ref{RHPPsib})--(\ref{zerorep}), the relation 
between the associated Toeplitz determinant and the $\tau$-function is given by
\be\label{tau4}
D_n(f(z)) = D_n(z_1,\dots, z_m|\alpha, \beta) =
\tau(z_1,\dots, z_m|\alpha, \beta) \prod_{j<k}(z_k-z_j)^{-2\alpha_j\alpha_k}
\prod_{j=0}^{m}z_j^{-2\alpha_j\lb}.
%\frac{\partial \ln D_n}{\partial z_k} = \frac{\partial \ln\tau}{\partial z_k} 
%+ 2\sum_{j=0, j\neq k}\frac{\alpha_k\alpha_j}{z_j-z_k}
%-2\kappa\frac{\alpha_k}{z_k}.
\ee

A direct substitution of formulae (\ref{Psij0}), (\ref{Psij01}), and (\ref{zerorep}) into
equations (\ref{fuchsian}) and (\ref{fuchs3}) yields
\begin{align}\label{tau010}
 \mbox{trace}\,
A_k\frac{\partial\wt{\Phi}_k(z_k)}{\partial z}
\wt{\Phi}^{-1}_k(z_k) 
&= \sum_{j=0\atop j\neq k}^{m}\frac{\mbox{trace}\,A_jA_k}{z_k-z_j} +\frac{1}{z_k}\mbox{trace}\,BA_k\\
&=  \sum_{j=0}^{m}\mbox{trace}\,
A_j\frac{\partial\wt{\Phi}_j(z_j)}{\partial z_k}
\wt{\Phi}^{-1}_j(z_j) + \mbox{trace}\,
B\frac{\partial\wt{\Phi}^{(0)}(0)}{\partial z_k}
\Bigl(\wt{\Phi}^{(0)}(0)\Bigr)^{-1}.\nonumber
\end{align}
Hence equation (\ref{tau3}) can be written as follows:
\begin{lemma} Let the Riemann-Hilbert problem for $\Phi$ be solvable.
For any $k=0,1,\dots,m$,
\be\la{tau55}
\frac{\partial \ln D_n}{\partial z_k} = \sum_{j=0}^{m}\mathrm{trace}\,
A_j\frac{\partial\wt{\Phi}_j(z_j)}{\partial z_k}
\wt{\Phi}^{-1}_j(z_j) + \mathrm{trace}\,
B\frac{\partial\wt{\Phi}^{(0)}(0)}{\partial z_k}
\Bigl(\wt{\Phi}^{(0)}(0)\Bigr)^{-1}
+ 2\sum_{j=0\atop j\neq k}^m\frac{\alpha_k\alpha_j}{z_j-z_k}
-2\lb\frac{\alpha_k}{z_k}.
\ee
\end{lemma}
What is of interest here is that the differential identities (\ref{diffid1}), which play a very
important role in the main text,  can be written in a
matrix form closely related to (\ref{tau55}). For simplicity, we present only
the case $\al_j\neq 0$. We have
\begin{lemma} Let the Riemann-Hilbert problem for $\Phi$ be solvable.
Let $\al_j\neq 0$, $j=0,\dots,m$. Then for any $k=0,1,\dots,m$,
\begin{multline}\la{tau5}
\frac{\partial \ln D_n}{\partial \alpha_k} = \sum_{j=0}^{m}\mathrm{trace}\,
A_j\frac{\partial\wt{\Phi}_j(z_j)}{\partial \alpha_k}
\wt{\Phi}^{-1}_j(z_j) + \mathrm{trace}\,
B\frac{\partial\wt{\Phi}^{(0)}(0)}{\partial \alpha_k}
\Bigl(\wt{\Phi}^{(0)}(0)\Bigr)^{-1}\\
 -2\sum_{j=0\atop j\neq k}^m\alpha_j\ln(z_j-z_k)
-2\lb\ln(-z_k) + \sum_{j=0}^{m}\alpha_j\ln z_j-n\ln z_k,
\end{multline}
\be\la{tau6}
\frac{\partial \ln D_n}{\partial \beta_k} =\sum_{j=0}^{m} \mathrm{trace}\,
A_j\frac{\partial\wt{\Phi}_j(z_j)}{\partial \beta_k}
\wt{\Phi}^{-1}_j(z_j) + \mathrm{trace}\,
B\frac{\partial\wt{\Phi}^{(0)}(0)}{\partial \beta_k}
\Bigl(\wt{\Phi}^{(0)}(0)\Bigr)^{-1}
 - \sum_{j=0}^{m}\alpha_j\ln z_j-n\ln z_k.
\ee
\end{lemma}
\begin{remark}
The significance of these equations is that they complement the
isomonodromy deformation formula (\ref{tau55}) by formulae
which describe the {\it monodromy} deformations  of the
$\tau$-function (represented by the Toeplitz determinants $D_n$).
(Equations (\ref{tau5}) and (\ref{tau6}) should be compared with
the general constructions of the recent paper \cite{BER2}.)
\end{remark}

\begin{proof}
Although straightforward, it is rather tedious to derive (\ref{tau5}) and (\ref{tau6}) 
from the differential identities (\ref{diffid1}). There is, however, an alternative
way to obtain (\ref{tau5}) and (\ref{tau6}) based on the direct analysis of the
Riemann-Hilbert problem (\ref{RHPYb})--(\ref{RHPYc}). 
First, we observe that the basic deformation formula for the Toeplitz 
determinant $D_n$, i.e. equations (\ref{premid}), (\ref{mid}) can be written,
using Eq. (2.4) in \cite{DIKt1}, in the form
\be\label{midnew}
{\partial \over\partial \ga}\ln D_n=
\frac{1}{2\pi i}\int_{C}z^{-n} 
\left(Y_{11}(z){dY_{21}(z)\over dz}
-{dY_{11}(z)\over dz}Y_{21}(z)\right)\frac{\partial f(z)}{\partial\ga} dz.
\ee
Here, as in Section \ref{diffid}, $\ga$ is either $\alpha_k$ 
of $\beta_k$. Second, by $\ga$--differentiating  the Riemann-Hilbert problem (\ref{RHPYb})--(\ref{RHPYc})
we easily obtain the following representation for the logarithmic derivative
$\frac{\partial Y(z)}{\partial\ga}Y^{-1}(z)$ of its solution (cf. Lemma 2.1 of \cite{BER2}):
\begin{align}\la{defgamma1}
X(z)\equiv\frac{\partial Y(z)}{\partial\ga}Y^{-1}(z)
&=\frac{1}{2\pi i}\int_{C}Y_{-}(z')\begin{pmatrix}0&1\cr
0&0\end{pmatrix}Y^{-1}_{+}(z')\frac{\partial f(z')}{\partial\ga}
(z')^{-n}\frac{dz'}{z'-z}\\
&= \frac{1}{2\pi i}\int_{C}\begin{pmatrix}-Y_{11}(z')Y_{21}(z')&Y^{2}_{11}(z')\cr
-Y^{2}_{21}(z')&Y_{11}(z')Y_{21}(z')\end{pmatrix}\frac{\partial f(z')}{\partial\ga}
(z')^{-n}\frac{dz'}{z'-z}.\nonumber
\end{align}
Now, from (\ref{Psitau}) and (\ref{fuchsian}) we see that
\begin{align}\nonumber
\frac{dY_{11}(z)}{dz} &= A_{11}(z)Y_{11}(z) +\Lambda^{-2}_{11} A_{12}(z)Y_{21}(z) - c(z)Y_{11}(z),\\
\frac{dY_{21}(z)}{dz} &= \Lambda^{2}_{11}A_{21}(z)Y_{11}(z) + A_{22}(z)Y_{21}(z) - c(z)Y_{21}(z),\nonumber
\end{align}
where $c(z) = \sum_{k=0}^{m}\frac{\alpha_k}{z-z_k} + \frac{\lb}{z}$.
This allows us to re-write (\ref{midnew}) as follows:
\begin{multline}\nonumber
{\partial \over\partial \ga}\ln D_n= \frac{1}{2\pi i}\int_{C}
z^{-n}\Bigl(Y^{2}_{11}(z)\Lambda^{2}_{11}A_{21}(z)  - \Lambda^{-2}_{11}A_{12}(z)Y^{2}_{21}(z)\\
+Y_{11}(z)Y_{21}(z)\left(A_{22}(z) - A_{11}(z)\right)\Bigr)\frac{\partial f(z)}{\partial\ga} dz\\
=\Lambda^{2}_{11} \sum_{k=0}^{m}A_{k,21}
\int_{C}Y^{2}_{11}(z)\frac{\partial f(z)}{\partial\ga} \frac{z^{-n}dz}{2\pi i(z-z_k)}
-\Lambda^{-2}_{11}\sum_{k=0}^{m}A_{k,12}
\int_{C}Y^{2}_{21}(z)\frac{\partial f(z)}{\partial\ga} \frac{z^{-n}dz}{2\pi i(z-z_k)}\\
+\sum_{k=0}^{m}\Bigl(A_{k,22}- A_{k,11}\Bigr)
\int_{C}Y_{11}(z)Y_{21}(z)\frac{\partial f(z)}{\partial\ga} \frac{z^{-n}dz}{2\pi i(z-z_k)}
+ \Lambda^{2}_{11}B_{21}\int_{C}Y^{2}_{11}(z)\frac{\partial f(z)}{\partial\ga} \frac{z^{-n}dz}{2\pi iz}\\
-\Lambda^{-2}_{11}B_{12}\int_{C}Y^{2}_{21}(z)\frac{\partial f(z)}{\partial\ga} \frac{z^{-n}dz}{2\pi iz}
+\Bigl(B_{22}- B_{11}\Bigr)
\int_{C}Y_{11}(z)Y_{21}(z)\frac{\partial f(z)}{\partial\ga} \frac{z^{-n}dz}{2\pi iz}.
\end{multline}
A comparison with (\ref{defgamma1}) yields the following {\it local}
representation for the $\ga$-derivative of $\ln D_n$:
\begin{multline}
{\partial \over\partial \ga}\ln D_n
=\Lambda^{2}_{11}\sum_{k=0}^{m}A_{k,21}X_{12}(z_k) + \Lambda^{-2}_{11}\sum_{k=0}^{m}A_{k,12}X_{21}(z_k)
+\sum_{k=0}^{m}A_{k,11}X_{11}(z_k)  + \sum_{k=0}^{m}A_{k,22}X_{22}(z_k)\\
+\Lambda^{2}_{11}B_{21}X_{12}(0) + \Lambda^{-2}_{11}B_{12}X_{21}(0)
+B_{11}X_{11}(0) + B_{22}X_{22}(0).
\end{multline}
%where we have  taken into account that all the matrices $A_k$ and $B$ are traceless.  
The last formula can be also written in the compact matrix form,
\be\label{defgamma2}
{\partial \over\partial \ga}\ln D_n = \sum_{k=0}^{m}\mbox{trace}\,\Lambda^{-1} A_k\Lambda X(z_k) +
\mbox{trace}\,\Lambda^{-1}B\Lambda X(0).
\ee
By evaluating (\ref{defgamma2}), with the help of the representation (\ref{Psij0}), 
one arrives at the  formulae (\ref{tau5}) and (\ref{tau6}). 
\end{proof}

\section*{Acknowledgements}
P. Deift was supported
in part by NSF grants \# DMS 0500923, \# DMS 1001886.
A. Its was supported
in part by NSF grant \#DMS-0701768 and
EPSRC grant \#EP/F014198/1. I. Krasovsky was
supported in part by EPSRC grants \#EP/E022928/1
and \#EP/F014198/1.

%%%%%%%%%%%%%%%%%%%%%%%%%%%%%%%%%%%%%%%%%%%%%%%
%%%%%%%%%%%%%%%%%%%%%%%%%%%%%%%%%%%%%%%%%%%%%%%%
%%%%%%%%%%%%%%%%%%%%%%%%%%%%%%%%%%%%%%%%%%%%%%%
%%%%%%%%%%%%%%%%%%%%%%%%%%%%%%%%%%%%%%%%%%%%%%

%%%%%%%%%%%%%%%%%%%%%%%%%%%%%%%%%%%%%%%
%%%%%%%%%%%%%%%%%%%%%%%%%%%%%%%%%%%%%%%
%%%%%%%%%%%%%%%%%%%%%%%%%%%%%%%%%%%%%%%
%%%%%%%%%%%%%%%%%%%%%%%%%%%%%%%%%%%%%%%
%%%%%%%%%%%%%%%%%%%%%%%%%%%%%%%%%%%%%%%
%%%%%%%%%%%%%%%%%%%%%%%%%%%%%%%%%%%%%%%
%%%%%%%%%%%%%%%%%%%%%%%%%%%%%%%%%%%%%%%
%%%%%%%%%%%%%%%%%%%%%%%%%%%%%%%%%%%%%%%
%%%%%%%%%%%%%%%%%%%%%%%%%%%%%%%%%%%%%%%
%%%%%%%%%%%%%%%%%%%%%%%%%%%%%%%%%%%%%%%
%%%%%%%%%%%%%%%%%%%%%%%%%%%%%%%%%%%%%%%

\end{document}